\documentclass[reqno]{amsart}
\usepackage{amssymb}
\usepackage{verbatim}

\begin{document}

\title[(Anti)symmetric multivariate trigonometric
functions]
{(Anti)symmetric multivariate trigonometric
functions and corresponding Fourier transforms}




{\Large

\noindent{\bf (Anti)symmetric multivariate trigonometric\\
functions and corresponding Fourier transforms}}

\bigskip

{\sf A.~Klimyk}

{\it Bogolyubov Institute for Theoretical Physics,
         Kiev 03680, Ukraine}

email: aklimyk@bitp.kiev.ua

\bigskip

{\sf J.~Patera}

{\it Centre de Recherches Math\'ematiques,
         Universit\'e de Montr\'eal,

         C.P.6128-Centre ville,
         Montr\'eal, H3C\,3J7, Qu\'ebec, Canada}

email: patera@crm.umontreal.ca

\bigskip

\noindent
{\it Abstract}.
Four families of special functions, depending on $n$ variables, are
studied. We call them symmetric and antisymmetric multivariate sine and
cosine functions. They are given as determinants or antideterminants of
matrices, whose matrix elements are sine or cosine functions of one
variable each. These functions are eigenfunctions of the Laplace
operator, satisfying specific conditions at the boundary of a certain
domain $F$ of the $n$-dimensional Euclidean space.
Discrete and continuous orthogonality on $F$ of the functions within
each family, allows one to introduce symmetrized and antisymmetrized
multivariate Fourier-like transforms, involving the symmetric and
antisymmetric multivariate sine and cosine functions.


\section{Introduction}
In mathematical and theoretical physics, very often we deal with
functions on the Euclidean space $E_n$ which are symmetric or
antisymmetric with respect to the permutation (symmetric) group
$S_n$. For example, such functions describe collections of identical
particles. Symmetric or antisymmetric solutions appear in the theory of
integrable systems. Characters of finite dimensional representations of
semisimple Lie algebras are symmetric functions. Moreover, according to
the Weyl formula for these characters, each such character is ratio of
antisymmetric functions.

The aim of this paper is to introduce, to describe and to study
symmetrized and antisymmetrized multivariate sine and cosine
functions and the corresponding Fourier transforms. Antisymmetric
multivariate sine and cosine functions (we denote them by ${\rm
sin}^-_\lambda(x)$, ${\rm cos}^-_\lambda(x)$,
$\lambda=(\lambda_1,\lambda_2, \dots,\lambda_n)$,
$x=(x_1,x_2,\dots,x_n)$) are determinants of $n\times n$ matrices,
whose entries are sine or cosine functions of one variable,
$$
{\rm sin}^-_\lambda(x)= \det (\sin 2\pi\lambda_ix_j)_{i,j=1}^n\,,
\qquad {\rm cos}^-_\lambda(x)= \det(\cos 2\pi
\lambda_ix_j)_{i,j=1}^n\,.
$$
These functions are antisymmetric in variables
$x_1,x_2,\dots,x_n$ and in parameters $\lambda_1,\lambda_2, \dots,\lambda_n$.
These antisymmetricities follow from antisymmetricity of
a determinant of a matrix under permutation of rows or
of columns.

  Symmetric multivariate sine and cosine functions ${\rm
sin}^+_\lambda(x)$, ${\rm cos}^+_\lambda(x)$ are antideterminants of
the same $n\times n$ matrices (a definition of antideterminants see
below). These functions are symmetric in variables and in
parameters. This symmetricity follows from symmetricity of an
antideterminant of a matrix under permutation of rows or of columns.

As in the case of sine and cosine functions of one variable, we
may consider three types of symmetric and antisymmetric multivariate
trigonometric functions:
  \medskip

(a) functions ${\rm sin}^\pm_m(x)$ and  ${\rm cos}^\pm_m(x)$ such
that $m=(m_1,m_2, \dots,m_n)$, $m_i\in \mathbb{Z}$; these functions
determine series expansions in multivariate sine and cosine
functions;
  \medskip

(b) functions ${\rm sin}^\pm_\lambda(x)$ and ${\rm
cos}^\pm_\lambda(x)$ such that $\lambda=(\lambda_1,\lambda_2,
\dots,\lambda_n)$, $\lambda_i\in \mathbb{R}$; these functions
determine sine and cosine integral Fourier transforms;
  \medskip

(c) functions ${\rm sin}^\pm_\lambda(x)$ and ${\rm
cos}^\pm_\lambda(x)$, where $x=(x_1,x_2, \dots,x_n)$ takes a finite
set of values; these functions determine multivariate finite sine
and cosine Fourier transforms.
  \medskip

  Functions (b) are antisymmetric (symmetric) with respect to elements
of the permutation group $S_n$. (Anti)symmetries of functions (a) are
described by a wider group, since sine and cosine functions $\sin 2\pi mx$, $\cos 2\pi
mx$, $m\in \mathbb{Z}$, are invariant with respect to shifts $x\to
x+k$, $k\in \mathbb{Z}$. (Anti)symmetries of functions (a)
are described by elements of the extended affine symmetric group
$\tilde S_n^{\rm aff}$ which is a product of the groups $S_n$, $T_n$
and $Z^n_2$, where $T_n$ consists of shifts of $E_n$ by vectors
$r=(r_1,r_2,\dots,r_n)$, $r_j\in  \mathbb{Z}$, and $Z^n_2$ is a product
of $n$ copies of the group $Z_2$ of order 2 generated by the reflection
with respect to the zero point. A fundamental domain $F(\tilde S_n^{\rm
aff})$ of the group $\tilde S_n^{\rm aff}$ is a certain bounded set in
$\mathbb{R}^n$.

The functions ${\rm cos}^+_m(x)$ with $m=(m_1,m_2, \dots,m_n)$,
$m_i\in\mathbb{Z}$, give solutions of the Neumann boundary value
problem on a closure of the fundamental domain $F(\tilde S_n^{\rm
aff})$. The functions ${\rm sin}^-_m(x)$ with $m=(m_1,m_2,
\dots,m_n)$, $m_i\in \mathbb{Z}$, are solutions of the Laplace
equation $\Delta f=\mu f$ on the domain $F(\tilde S_n^{\rm aff})$
vanishing on the boundary $\partial F(\tilde S_n^{\rm aff})$ of
$F(\tilde S_n^{\rm aff})$ (Dirichlet boundary problem). The
functions (b) are also solutions of the Laplace equation.

Functions on the fundamental domain $F(\tilde S_n^{\rm aff})$ can be
expanded into series in the functions (a). These expansions are an
analogue of the usual sine and cosine Fourier series for functions of
one variable. Functions (b) determine antisymmetrized and symmetrized
sine and cosine Fourier integral transforms on the fundamental domain
$F(\tilde S_n)$ of the extended symmetric group $\tilde S_n=S_n\times Z_2^n$.
This domain consists of points $x\in E_n$ such that $x_1>x_2>\cdots
 >x_n>0$.

Functions (c) are used to determine antisymmetric or symmetric finite
(that is, on a finite set) trigonometric multivariate Fourier
transforms. These finite trigonometric transforms are given on grids
consisting of points of the fundamental domain $F(\tilde S_n^{\rm
aff})$.

Symmetric and antisymmetric sine and cosine functions, which
are studied in this paper, are closely related to symmetric and
antisymmetric orbit functions defined in Refs. 1 and 2 and studied in
detail in Refs. 3 and 4. They are connected with orbit
functions corresponding to the Dynkin--Coxeter diagrams
of semisimple Lie algebras of rank $n$.
Discrete orbit function transforms, corresponding to
Dynkin--Coxeter diagrams of low rank, were studied and exploited in
rather useful applications${}^{5-13}$. Clearly, our multivariate sine
and cosine transforms can be applied under solution of the same problems,
that is, of the problems formulated on grids or lattices.

The exposition of the theory of orbit functions in Refs. 3 and 4
strongly depends on the theory of Weyl groups, properties of root
systems, etc. In this paper we avoid this dependence (moreover, some of
our functions cannot be treated by using Weyl gorups and root systems).
We use only the permutation (symmetric) group, its extension, and
properties of determinants and antideterminants.

It is well-known that
the determinant $\det (a_{ij})_{i,j=1}^n$ of an $n\times n$ matrix
$(a_{ij})_{i,j=1}^n$ is defined as
  \begin{align}\label{det}
\det (a_{ij})_{i,j=1}^n= & \sum_{w\in S_n} (\det w)
a_{1,w(1)}a_{2,w(2)}\cdots a_{n,w(n)}
\notag \\
= &\sum_{w\in S_n} (\det w)
a_{w(1),1}a_{w(2),2}\cdots a_{w(n),n},
\end{align}
where $S_n$ is the permutation (symmetric) group of $n$ symbols
$1,2,\dots,n$, the set $(w(1),w(2),\dots$, $w(n))$ means the set
$w(1,2,\dots,n)$, and  $\det w$ denotes a determinant of the transform
$w$, that is, $\det w=1$ if $w$ is an even permutation and $\det w=-1$
otherwise. Along with a determinant, we use an antideterminant $\det^+$
of the matrix $(a_{ij})_{i,j=1}^n$ which is defined as a sum of
all summands entering into the expression for a determinant,
taken with the sign +,
  \begin{align}\label{dett}
{\det}^+ (a_{ij})_{i,j=1}^n= & \sum_{w\in S_n}
a_{1,w(1)}a_{2,w(2)}\cdots
a_{n,w(n)}
\notag \\
= & \sum_{w\in S_n} a_{w(1),1}a_{w(2),2}\cdots
a_{w(n),2}.
\end{align}

Symmetrized or antisymmetrized multivariate sine and cosine functions
were mentioned in Refs. 14 and 15. In this paper we investigate in detail
these multivariate functions and derive the
corresponding continuous and finite Fourier transforms. Note that
in Ref. 16 we have studied symmetric and antisymmetric exponential functions.

\section{Symmetric and antisymmetric multivariate sine and
cosine functions}

Antisymmetric multivariate sine functions ${\rm sin}^-_{\lambda}(x)$
on ${\mathbb{R}}^n$ are determined by $n$ real numbers
$\lambda\equiv (\lambda_1,\lambda_2,\dots,\lambda_n)$ and are given
by the formula
 \begin{align}\label{det-B}
{\rm sin}^-_{\lambda}(x)\equiv & \; {\rm
sin}^-_{(\lambda_1,\lambda_2,\dots ,\lambda_{n})}(x_1,x_2,\dots,x_n)
:= \det \left(
 \sin 2\pi \lambda_ix_j  \right)_{i,j=1}^{n}
\notag \\
=&\;  \det \left(
 \begin{array}{cccc}
 \sin 2\pi\lambda_1x_1& \sin 2\pi\lambda_1x_2&\cdots & \sin 2\pi\lambda_1x_{n}\\
  \sin 2\pi\lambda_2x_1& \sin 2\pi\lambda_2x_2&\cdots & \sin 2\pi\lambda_2x_{n}\\
  \cdots & \cdots &
\cdots & \cdots \\
 \sin 2\pi\lambda_{n}x_1& \sin 2\pi\lambda_{n}x_2&\cdots & \sin 2\pi\lambda_{n}x_{n}
 \end{array} \right)
\notag\\
\equiv & \sum_{w\in S_n} (\det w)
 \sin 2\pi\lambda_1x_{w(1)} \sin 2\pi\lambda_2x_{w(2)}\cdots
 \sin 2\pi\lambda_nx_{w(n)}
\notag\\
= & \sum_{w\in S_n} (\det w)
 \sin 2\pi\lambda_{w(1)}x_{1} \sin 2\pi\lambda_{w(2)}x_{2}\cdots
 \sin 2\pi\lambda_{w(n)}x_{n},
 \end{align}
where $(w(1),w(2),\dots$, $w(n))$ means the set
$w(1,2,\dots,n)$ and  $\det w$ denotes a determinant of the transform $w$,
$\det w=\pm 1$.

A special case of the antisymmetric multivariate sine functions is
when $\lambda_i$ are integers; in this case we write
$(m_1,m_2,\dots, m_n)$ instead of $(\lambda_1,\lambda_2,\dots,m_n)$,
\[
{\rm sin}^-_{(m_{1},m_2,\dots ,m_{n})}(x) =  \det \left(
 \sin 2\pi m_ix_j  \right)_{i,j=1}^{n},\ \ \ \ m\in \mathbb{Z}.
\]

If $n=2$ we have
\begin{align}  
{\rm sin}^-_{(\lambda_1,\lambda_2)}(x_1,x_2) = &
 \sin \pi\lambda_1x_1\;  \sin \pi\lambda_2x_2-
 \sin \pi\lambda_1x_2\;  \sin \pi\lambda_2x_1 \notag\\
= &\; \frac12 \cos 2\pi(\lambda_1x_1 - \lambda_2x_2)- \frac12 \cos
2\pi(\lambda_1x_1 +\lambda_2x_2)
 \notag\\
& \; - \frac12 \cos 2\pi(\lambda_1x_2  -\lambda_2x_1)+ \frac12 \cos
2\pi(\lambda_1x_2 + \lambda_2x_1).
\notag
 \end{align}

Antisymmetric multivariate cosine functions on ${\mathbb{R}}^n$ are
determined by $n$ real numbers $\lambda=(\lambda_1,\lambda_2,\dots,
\lambda_n)$ and are given by the formula
 \begin{align}\label{det-D-2}
{\rm cos}^-_{\lambda}(x)\equiv & \; {\rm
cos}^-_{(\lambda_{1},\lambda_2,\dots
,\lambda_{n})}(x_1,x_2,\dots,x_n)
 :=  \det \left(  \cos 2\pi \lambda_ix_j  \right)_{i,j=1}^{n}
\notag\\
\equiv & \sum_{w\in S_n} (\det w)
 \cos 2\pi\lambda_1x_{w(1)} \cos 2\pi\lambda_1x_{w(2)}\cdots
 \cos 2\pi\lambda_1x_{w(n)}
\notag\\
= & \sum_{w\in S_n}   (\det w)
 \cos 2\pi\lambda_{w(1)}x_{1} \cos 2\pi\lambda_{w(2)}x_{2}\cdots
 \cos 2\pi\lambda_{w(n)}x_{n}.
 \end{align}

As in the case of sine functions, we separate the special case of
the antisymmetric multivariate cosine functions when $\lambda_i$ are
integers and write $(m_1,m_2,\dots, m_n)$ instead of
$(\lambda_1,\lambda_2,\dots,\lambda_n)$.

If $n=2$, we get
\begin{align}  
{\rm cos}^-_{(\lambda_1,\lambda_2)}(x_1,x_2) =& \; \frac12 \cos
2\pi(\lambda_1x_1 - \lambda_2x_2)+ \frac12 \cos 2\pi(\lambda_1x_1
+\lambda_2x_2)
 \notag\\
 & \; -  \frac12 \cos 2\pi(\lambda_1x_2  -\lambda_2x_1)- \frac12 \cos
2\pi(\lambda_1x_2 + \lambda_2x_1).
\notag
 \end{align}

Symmetric multivariate sine functions on ${\mathbb{R}}^n$ are
determined by $n$ real numbers $(\lambda_1,\lambda_2,\dots,\lambda_n)$ and
are given by
 \begin{align}\label{det-B-S}
{\rm sin}^+_{\lambda}(x) \equiv & \; {\rm
sin}^+_{(\lambda_{1},\lambda_2,\dots
,\lambda_{n})}(x_1,x_2,\dots,x_n) :=  {\det}^+ \left(
 \sin 2\pi \lambda_ix_j  \right)_{i,j=1}^{n}
\notag\\
 =&\;
\sum_{w\in S_n} \sin 2\pi m_{1}x_{w(1)} \cdot {\dots}\cdot \sin 2\pi
m_{n}x_{w(n)}
\notag\\
 =&\;
\sum_{w\in S_n} \sin 2\pi m_{w(1)}x_1 \cdot {\dots}\cdot \sin 2\pi
m_{w(n)}x_n,
 \end{align}
If $n=2$ we have
\begin{align}
{\rm sin}^+_{(\lambda_1,\lambda_2)}(x_1,x_2) =& \frac12 \cos
2\pi(\lambda_1x_1 - \lambda_2x_2)- \frac12 \cos 2\pi(\lambda_1x_1
+\lambda_2x_2)
 \notag\\
& +  \frac12 \cos 2\pi(\lambda_1x_2  -\lambda_2x_1)- \frac12 \cos
2\pi(\lambda_1x_2 + \lambda_2x_1).
\notag
 \end{align}

Symmetric multivariate cosine functions on ${\mathbb{R}}^n$ are
determined by $n$ real numbers $\lambda=(\lambda_1,\lambda_2,\dots,\lambda_n)$ and
are given by the formula
 \begin{align}\label{det-B-C}
{\rm cos}^+_{\lambda}(x) =&\;  {\det}^+ \left(
 \cos 2\pi \lambda_ix_j  \right)_{i,j=1}^{n}\notag \\
 =&\;
\sum_{w\in S_n} \cos 2\pi m_{1}x_{w(1)} \cdot {\dots}\cdot \cos 2\pi
m_{n}x_{w(n)} \notag\\
 =&\;
\sum_{w\in S_n} \cos 2\pi m_{w(1)}x_1 \cdot {\dots}\cdot \cos 2\pi
m_{w(n)}x_n.
 \end{align}
If $n=2$, we get
\begin{align}
{\rm cos}^+_{(\lambda_1,\lambda_2)}(x_1,x_2) =& \frac12 \cos
2\pi(\lambda_1x_1 - \lambda_2x_2)+ \frac12 \cos 2\pi(\lambda_1x_1
+\lambda_2x_2)
 \notag\\
 & +  \frac12 \cos 2\pi(\lambda_1x_2  -\lambda_2x_1)+ \frac12 \cos
2\pi(\lambda_1x_2 + \lambda_2x_1).
\notag
 \end{align}

A special case of the symmetric multivariate sine and cosine
functions is when $\lambda_i$ are integers; in this case we write
$(m_1,m_2,\dots, m_n)$ instead of $(\lambda_1,\lambda_2,\dots,m_n)$.

\section{Extended affine symmetric group and fundamental domains}

It follows from the definitions of the functions ${\rm
sin}^{\pm}_\lambda(x)$ and ${\rm cos}^{\pm}_\lambda(x)$ that a
permutation of variables $x_1,x_2,\dots,x_n$ is equivalent to the
same permutation of the corresponding rows in the matrices from the
definition. Then due to properties of determinants and
antideterminants of matrices under permutations of rows, the
functions ${\rm sin}^\pm_\lambda(x)$ and ${\rm cos}^\pm_\lambda(x)$
are symmetric or antisymmetric with respect to the permutation group
$S_n$, that is,
 \[
{\rm sin}^+_\lambda(wx)={\rm sin}^+_\lambda(x), \ \ \ \ {\rm
cos}^+_\lambda(wx)={\rm cos}^+_\lambda(x),\ \ \ \ w\in S_n,
\]  \[
{\rm sin}^-_\lambda(wx)=(\det w){\rm sin}^-_\lambda(x), \ \ {\rm
cos}^-_\lambda(wx)=(\det w){\rm cos}^-_\lambda(x),\ \ w\in S_n.
\]

These functions admit additional (anti)invariances with respect to
changing signs of coordinates $x_1,x_2,\dots,x_n$. Let $\varepsilon_i$ denote
the operation of a change of a sign of the coordinate $x_i$.
Since $\sin 2\pi \varepsilon_ix_i\lambda_j=
-\sin 2\pi x_i\lambda_j$ and $\cos 2\pi \varepsilon_ix_i\lambda_j=
\cos 2\pi x_i\lambda_j$, we have
\begin{align}\label{change-1}
{\rm sin}^+_\lambda(\varepsilon_ix)=- {\rm sin}^+_\lambda(x),\ \ \ &
{\rm cos}^+_\lambda(\varepsilon_ix)=
{\rm cos}^+_\lambda(x),\\
{\rm sin}^-_\lambda(\varepsilon_ix)=- {\rm sin}^-_\lambda(x),\ \ \ &
{\rm cos}^-_\lambda(\varepsilon_ix)= {\rm cos}^-_\lambda(x).
\end{align}
We denote the group generated by changes of coordinate signs by
$Z_2^n$, where $Z_2$ is the group of changes of sign of one
coordinate.

The group $\tilde S_n=S_n\times Z_2^n$ is called the {\it extended
symmetric group}. It is a group of (anti)symmetries for the
functions ${\rm sin}^\pm_\lambda(x)$ and ${\rm cos}^\pm_\lambda(x)$.

We have the same (anti)symmetries under changes of signs in numbers
$\lambda_1,\lambda_2,\dots$, $\lambda_n$. In order to avoid these
(anti)symmetries, we may assume that all
coordinates $x_1,x_2,\dots,x_n$ and all numbers
$\lambda_1,\lambda_2,\dots,\lambda_n$ are non-negative.

The trigonometric functions ${\rm sin}^+_m$ and ${\rm cos}^+_m$,
determined by integral $m=(m_1$, $m_2,\dots, m_n)$, admit additional
symmetries related to periodicity of the sine and cosine functions
$\sin 2\pi r y$, $\cos 2\pi r y$ of one variable for integral values
of $r$. These symmetries are described by the discrete group of
shifts by vectors
\[
r_1{\bf e}_1+r_2{\bf e}_2+\cdots +r_n{\bf e}_n,\ \ \ r_i\in
\mathbb{Z},
\]
where ${\bf e}_1,{\bf e}_2,\dots, {\bf e}_n$ are the unit vectors in
directions of the corresponding axes. We denote this group by
$T_n$. Permutations of $S_n$, the operations $\varepsilon_i$ of
changes of coordinate signs, and shifts of $T_n$ generate a group
which is denoted as $\tilde S_n^{\rm aff}$ and is called the
{\it extended affine symmetric group}. (The group generated by
permutations of $S_n$ and by shifts of $T_n$ generate a group
which is denoted as $S_n^{\rm aff}$ and is called the
{\it affine symmetric group}). Thus, the group $\tilde S_n^{\rm aff}$
is a product of its subgroups,
\[
\tilde S_n^{\rm aff}=S_n\times Z_2^n\times T_n,
\]
where $T_n$ is an invariant subgroup, that is, $wtw^{-1}\in T_n$ and
$\varepsilon_i t\varepsilon^{-1}_i\in T_n$ for
$w\in S_n$, $\varepsilon_i\in Z_2$, $t\in T_n$.

An open connected simply connected set $F\subset \mathbb{R}^n$ is
called a {\it fundamental domain} for the group $\tilde S_n^{\rm aff}$ (for
the group $\tilde S_n$) if it does not contain equivalent points (that is,
points $x$ and $x'$ such that $x'=gx$, where $g$ is an element of
$\tilde S_n^{\rm aff}$ or $\tilde S_n$) and if its closure contains at
least one point from each $\tilde S_n^{\rm aff}$-orbit (from each
$\tilde S_n$-orbit). Recall that an $\tilde S_n^{\rm aff}$-orbit of a point $x\in
\mathbb{R}^n$ is the set of points $wx$, $w\in \tilde S_n^{\rm aff}$.

It is evident that the set $D_+^+$ of all points
$x=(x_1,x_2,\dots,x_n)$ such that
\begin{equation}\label{aff-dom}
x_1>x_2>\cdots >x_n>0
\end{equation}
is a fundamental domain for the group $\tilde S_n$ (we denote it as
$F(\tilde S_n)$). The set of points $x=(x_1,x_2,\dots,x_n)\in D^+_+$ such
that
\begin{equation}\label{aff-dom-}  {\textstyle
\frac 12 >x_1>x_2>\cdots >x_n>0    }
\end{equation}
is a fundamental domain for the extended affine
symmetric group $\tilde S^{\rm aff}_n$
(we denote it as $F(\tilde S^{\rm aff}_n)$).
 \medskip

{\it Remark.} It may be seemed that instead of $\frac 12$ in \eqref{aff-dom-}
there must be 1. However, in the group $\tilde S_n^{\rm aff}$ there exists
the reflection of $\mathbb{R}^n$ with respect to the hyperplane
$x_1=\frac 12$. This reflection coincides with the transform
\[
 x\to \varepsilon_1x+{\bf e}_1,
\]
where $\varepsilon_1$ is a change of a sign of the coordinate $x_1$ (note
that $\varepsilon_1$ is the reflection of $\mathbb{R}^n$ with respect to the hyperplane
$x_1=0$). The transform $x\to \varepsilon_1x+{\bf e}_1$ leaves coordinates
$x_2,x_3,\dots,x_n$ invariant and does not move the hyperplane $x_1=\frac 12$,
that is, it is a reflection. Therefore, the domain $1>x_1>x_2>\cdots
>x_n>0$ consists of two copies of the fundamental domain
$F(\tilde S^{\rm aff}_n)$.

In the group $\tilde S_n^{\rm aff}$ there also exist
the reflections of $\mathbb{R}^n$ with respect to the hyperplanes
$x_i=\frac 12$, $i=2,3,\dots,n$.
 \medskip

As we have seen, the symmetric multivariate trigonometric functions
${\rm sin}^+_\lambda(x)$ and ${\rm cos}^+_\lambda(x)$ are symmetric
with respect to the symmetric group $S_n$ and behave according to
formula \eqref{change-1} under a change of coordinate signs. This
means that it is enough to consider these functions only on the
closure of the fundamental domain $F(\tilde S_n)$, that is, on the
set $D_+$ of points $x$ such that
\[
 x_1\ge x_2\ge \cdots \ge x_n\ge 0.
\]
Values of these functions on other
points are received by using symmetricity.

Symmetricity of ${\rm sin}^+_m$ and ${\rm cos}^+_m$ with integral
$m=(m_1,m_2,\dots,m_n)$ with respect to the extended affine
symmetric group $\tilde S_n^{\rm aff}$,
\begin{equation}\label{aff-sym-aff}
{\rm sin}^+_m(wx+r)={\rm sin}^+_m(x), \ \ {\rm cos}^+_m(wx+r)={\rm
cos}^+_m(x)\ \ \ \
 w\in S_n,\ \ \ r\in T_n,
\end{equation}
\begin{equation}\label{aff-sym-aff-2}
{\rm sin}^+_m(\varepsilon_i x)=-{\rm sin}^+_m(x), \ \ {\rm
cos}^+_m(\varepsilon_i x)={\rm cos}^+_m(x),\ \ \ \ \varepsilon_i\in
Z_2,
\end{equation}
means that we may consider the functions ${\rm sin}^+_m$ and ${\rm
cos}^+_m$ on the closure of the fundamental domain $F(\tilde S^{\rm
aff}_n)$, that is, on the set of points
\[   {\textstyle
\frac 12\ge x_1\ge x_2\ge \cdots \ge x_n\ge 0.  }
\]
Values of these functions on other points are obtained by
using relations \eqref{aff-sym-aff} and \eqref{aff-sym-aff-2}.

The functions ${\rm sin}^-_\lambda(x)$ and ${\rm cos}^-_\lambda(x)$
are antisymmetric with respect to the extended symmetric group
$\tilde S_n$,
\[
{\rm sin}^-_\lambda(wx)=(\det w){\rm sin}^-_\lambda(x),\
 \ {\rm cos}^-_\lambda(wx)=(\det w){\rm cos}^-_\lambda(x),
\ \ \ \ w\in S_n,
\] \[
{\rm sin}^-_\lambda(\varepsilon_i x)=-{\rm sin}^-_\lambda(x), \ \
{\rm cos}^-_\lambda(\varepsilon_i x)={\rm cos}^-_\lambda(x),\ \ \ \
\varepsilon_i\in Z_2.
\]
For this reason, we may consider these functions only
on the fundamental domain $F(\tilde S_n)$.

The antisymmetric sine and cosine functions ${\rm sin}^-_m(x)$ and
${\rm cos}^-_m(x)$ with integral $m=(m_1,m_2,\dots, m_n)$ also admit
additional symmetries related to the periodicity of the usual sine
and cosine functions. These symmetries are described by elements of
the extended affine symmetric group $\tilde S^{\rm aff}_n$. For
$w\in S_n$, $r\in T_n$ and $\varepsilon_i\in Z_2$ we have
\begin{equation}\label{aff-sym-af}
{\rm sin}^-_m(wx+r)=(\det w){\rm sin}^-_m(x), \ \ \ {\rm
cos}^-_m(wx+r)=(\det w){\rm cos}^-_m(x),
\end{equation}
\begin{equation}\label{aff-sym-af-2}
{\rm sin}^-_m(\varepsilon_i x)=-{\rm sin}^-_m(x), \ \ \ {\rm
cos}^-_m(\varepsilon_i x)={\rm cos}^-_m(x),
\end{equation}
that is, it is enough
to consider these functions only on the
closure of the fundamental domain $F(\tilde S^{\rm aff}_n)$. Values of
these functions on other points are obtained by using
the relations \eqref{aff-sym-af} and \eqref{aff-sym-af-2}.

\section{Properties}
Symmetricity and antisymmetricity of symmetric and antisymmetric
multivariate trigonometric functions is a main property of these
functions. However, they possess many other interesting properties.
 \medskip

{\bf Behavior on boundary of fundamental domain.} The symmetric and
antisymmetric functions ${\rm sin}^{\pm}_\lambda(x)$ and ${\rm
cos}^{\pm}_\lambda(x)$ are finite sums of products of sine or cosine
functions. Therefore, they are continuous functions of
$x_1,x_2,\dots,x_n$ and have continuous derivatives of all orders on
$\mathbb{R}^n$.

The boundary $\partial F(\tilde S_n)$
of the fundamental domain $F(\tilde S_n)$ consists of points
of $F(\tilde S_n)$ which belong at least to one of the hyperplanes
\[
x_1=x_2,\ \ x_2=x_3,\ \ \cdots, \ \ x_{n-1}=x_n,\ \ x_n=0.
\]
The set of points of the boundary, determined by the hyperplane
$x_i=x_{i+1}$ or by the hyperplane $x_n=0$, is called a wall of the fundamental domain.

Since for $x_i=x_{i+1}$, $i=1,2,\dots,n-1$, the matrix $(\sin 2\pi
\lambda_ix_j)_{i,j=1}^n$ has two coinciding columns, then $\det
(\sin 2\pi \lambda_ix_j)_{i,j=1}^n=0$ in these cases. Clearly, we
also have $\det (\sin 2\pi \lambda_ix_j)_{i,j=1}^n=0$ for $x_n=0$.
Thus, {\it the function ${\rm sin}^-_\lambda(x)$ vanishes on the
boundary $\partial F(\tilde S_n)$}.

It is shown similarly that the function ${\rm cos}^-_\lambda(x)$
vanishes on the walls $x_i=x_{i+1}$, $i=1,2,\dots,n-1$, of the
boundary and $\partial {\rm cos}^-_\lambda(x)/\partial x_n$ vanishes
on the wall $x_n=0$.

The antideterminant ${\det}^+ (\cos 2\pi \lambda_ix_j)_{i,j=1}^n$
does not change under permutation of two coinciding columns which
appear on the walls $x_i=x_{i+1}$, $i=1,2,\dots,n-1$, of the
boundary $\partial F(\tilde S_n)$. Moreover, the function ${\rm
cos}^+_\lambda(x)$ is invariant under reflections $r_i$ with respect
to the hyperplanes $x_i=x_{i+1}$, $i=1,2,\dots,n-1$ (these
reflections lead to permutations of the corresponding coordinates
$x_i$ and $x_{i+1}$). We also have $\partial\, {\rm
cos}^+_\lambda(x)/\partial\, x_n=0$ for $x_n=0$. These assertions
means that
\[
\frac{\partial\, {\rm cos}^+_\lambda(x)}{\partial\, {\bf n}}
\left\vert_{\partial F(\tilde S_n)}=0, \right.
\]
where $\bf n$ is the normal to the boundary $\partial F(\tilde S_n)$.

For ${\rm sin}^+_\lambda(x)$ we have that $\partial\, {\rm
sin}^+_\lambda(x)/\partial\, {\bf n}$ vanishes on the walls
$x_i=x_{i+1}$, $i=1,2,\dots,n-1$, and ${\rm sin}^+_\lambda(x)=0$ on
the wall $x_n=0$.
 \medskip

{\bf Scaling symmetry.} Let $\lambda=(\lambda_1,\lambda_2,\dots,\lambda_n)$.
For $c\in \mathbb{R}$ we denote by $c\lambda$ the set $(c\lambda_1,c\lambda_2,\dots,
c\lambda_n)$. We have
\[
{\rm sin}^-_{c\lambda}(x)= \det(\sin 2\pi c\lambda_ix_j)_{i,j=1}^n=
 {\rm sin}^-_{\lambda}(cx).
\]
This equality expresses the {\it scaling symmetry of
antisymmetric sine functions.} It is shown similarly that
\[
{\rm cos}^-_{c\lambda}(x)={\rm cos}^-_{\lambda}(cx),\ \ \ {\rm
sin}^+_{c\lambda}(x)={\rm sin}^+_{\lambda}(cx),
 \ \ \
{\rm cos}^+_{c\lambda}(x)={\rm cos}^+_{\lambda}(cx).
\]

{\bf Duality.} It follows from formulas for antisymmetric
sine and cosine functions that they do not change under permutation
$(\lambda_1,\lambda_2,\dots,\lambda_n)\leftrightarrow
(x_1,x_2,\dots,x_n)$, $x_i\ne x_j$, $\lambda_i\ne \lambda_j$,
that is, we have
 \[
{\rm sin}^-_\lambda(x)= {\rm sin}^-_x(\lambda),\ \ \ \ {\rm
cos}^-_\lambda(x)= {\rm cos}^-_x(\lambda).
 \]
These relations expresses the {\it duality} of antisymmetric sine
and cosine functions.

The duality is also true for symmetric trigonometric functions,
\[
{\rm sin}^+_\lambda(x)= {\rm sin}^+_x(\lambda),\ \ \ \ {\rm
cos}^+_\lambda(x)= {\rm cos}^+_x(\lambda).
\]

{\bf Orthogonality on the fundamental domain $F(\tilde S^{\rm
aff}_n)$.} Antisymmetric multivariate sine functions ${\rm sin}^-_m$
with $m=(m_1,m_2,\dots,m_n)\in D^+_+$, $m_j\in \mathbb{Z}$, are
orthogonal on $F(\tilde S_n^{\rm aff})$ with respect to the
Euclidean measure. We have
 \begin{equation}\label{symmmm-cos}
2^{2n} \int_{\overline{F(\tilde S^{\rm aff}_n)}} {\rm sin}^-_m
(x){\rm sin}^-_{m'}(x)dx= \delta_{m, m'} ,
 \end{equation}
where the closure ${\overline{F(\tilde S^{\rm aff}_n)}}$
of $F(\tilde S^{\rm aff}_n)$ consists of
points $x=(x_1,x_2, \dots,x_n)\in E_n$ such that
\[
 \frac12 \ge x_1\ge x_2\ge\cdots \ge x_n\ge 0.
\]
This relation follows from  orthogonality of the sine functions
$\sin 2\pi m_ix_j$ (entering into the definition of the function
${\rm sin}^-_m(x)$). Indeed, we have
\[
2^2 \int^{1/2}_0 \sin (2\pi kt)\, \sin (2\pi k't)\, dt =\delta_{kk'}.
\]
Therefore, if ${\sf T}$ is the set $[0,\frac12 ]^n$, then
\[
2^{2n} \int_{\sf T} {\rm sin}^-_m (x){\rm sin}^-_{m'}(x)dx= \vert
S_n\vert \delta_{m, m'} ,
\]
where $\vert S_n\vert$ is an order of the permutation group. Since
$F(\tilde S^{\rm aff}_n)$ covers the set ${\sf T}$ exactly
$\vert S_n\vert$ times, the formula \eqref{symmmm-cos} follows.

A similar orthogonality relation can be written down for the
antisymmetric multivariate cosine functions,
 \begin{equation}\label{symmmm-cos-2}
2^{2n} \int_{\overline{F(\tilde S^{\rm aff}_n)}} {\rm cos}^-_m(x)
{\rm cos}^-_{m'}(x)\, dx= \delta_{m, m'} .
 \end{equation}

For symmetric multivariate sine and cosine function we have the
orthogonality relations
\begin{equation}\label{ortog-ant-cos-2}
2^{2n}\int_{\overline{F(\tilde S^{\rm aff}_n)}} {\rm sin}^+_m(x)
{\rm sin}^+_{m'}(x)\, dx = |G_m|\delta_{m,m'} ,
 \end{equation}
\begin{equation}\label{ortog-ant-cos}
2^{2n}\int_{\overline{F(\tilde S^{\rm aff}_n)}} {\rm cos}^+_m(x)
{\rm cos}^+_{m'}(x)\, dx = |G_m|\delta_{m, m'} ,
 \end{equation}
where $m$ and $m'$ are such that $m_1\ge m_2\ge
\cdots\ge m_n\ge 0$, $m'_1\ge m'_2\ge \cdots\ge m'_n\ge 0$, $m_i,m'_j\in
\mathbb{Z}$, and $|G_m|$
is an order of the subgroup $G_m\subset S_n$ consisting of
elements leaving $m$ invariant.

{\bf Orthogonality of symmetric and antisymmetric
trigonometric functions.}
Let $w_i$ ($i=1,2,\dots,n-1$) be the permutation of coordinates $x_i$ and
$x_{i+1}$. We create the domain $F^{\rm ext}(\tilde S_n^{\rm aff})
=F(\tilde S_n^{\rm aff})\cup w_i
F(\tilde S_n^{\rm aff})$,
where $F(\tilde S_n^{\rm aff})$ is the
fundamental domain for the group $\tilde S_n^{\rm aff}$. Let $F^{\rm ext}$
be a closure of the domain $F^{\rm ext}(\tilde S_n^{\rm aff})$.
If $i=1$, then $F^{\rm ext}$ consists of points $x\in E_n$ such that
\[
\frac12 \ge x_1\ge x_2\ge\cdots \ge x_n\ge 0\ \ \ {\rm or}\ \ \
\frac12 \ge x_2\ge x_1\ge x_3\ge x_4\ge\cdots \ge x_n\ge 0.
\]

 Since for $m=(m_1,m_2,\dots,m_n)\in
 D_+\equiv \overline{D^+_+}$, $m_j\in \mathbb{Z}$, we have
\[
{\rm cos}^+_m(w_i x)={\rm cos}^+_m(x),\qquad {\rm sin}^-_m(w_i
x)=-{\rm sin}^-_m(x) ,
\]
then
 \begin{equation}\label{ortog-ant}
 \int_{F^{\rm ext}} {\rm sin}^-_m(x)
{\rm cos}^+_{m'}(x) dx=0.
 \end{equation}
Indeed, due to symmetry and antisymmetry of symmetric and
antisymmetric trigonometric functions, respectively, we have
\[
\int_{F^{\rm ext}} {\rm sin}^-_m(x) {\rm cos}^+_{m'}(x)dx
\] \[  \qquad
= \int_{\overline{F(\tilde S^{\rm aff}_n)}} {\rm sin}^-_m(x) {\rm
cos}^+_{m'}(x)dx + \int_{w_{i}{\overline{F(\tilde S^{\rm aff}_n)}}}
{\rm sin}^-_m(x) {\rm cos}^+_{m'}(x)dx
\] \[   \qquad
= \int_{\overline{F(\tilde S^{\rm aff}_n)}} {\rm sin}^-_m(x) {\rm
cos}^+_{m'}(x)dx+ \int_{\overline{F(\tilde S^{\rm aff}_n)}} (-{\rm
sin}^-_m(x)) {\rm cos}^+_{m'}(x)dx=0.
\]
For $n=1$ the orthogonality \eqref{ortog-ant} means the
orthogonality of the functions sine and cosine on the interval
$(0,2\pi)$.

The relation
 \begin{equation}\label{ortog-ant-3}
 \int_{F^{\rm ext}} {\rm sin}^+_m(x)
{\rm cos}^-_{m'}(x) dx=0
 \end{equation}
is proved similarly.

\section{Special cases}

For special values of
$\lambda=(\lambda_1,\lambda_2,\dots,\lambda_n)$ the function ${\rm
sin}^-_\lambda(x)$ can be represented in a form of product of
trigonometric functions of one variables. For $\lambda\equiv
\rho_1=(n,n-1,\dots, 1)$ we get
\begin{equation}\label{anttttt-3}
{\rm sin}^-_{\rho_1}(x)=\prod_{1\le i<j\le n} \sin \, \pi (x_i-x_j)
\sin \, \pi (x_i+x_j)
 \prod_{1\le i\le n} \sin 2\pi x_i.
\end{equation}
In order to prove this formula we have to represent the sine
functions of one variable in \eqref{det-B} and \eqref{anttttt-3} in
terms of exponential functions. Then we fulfil multiplications of
all expressions in \eqref{anttttt-3} and obtain ${\rm
sin}^-_{\rho_1}(x)$ in the form of sum of products of exponential
functions. Comparing this form with the expression \eqref{det-B} for
${\rm sin}^-_{\rho_1}(x)$ in terms of exponential functions we show
that formula \eqref{anttttt-3} is true.

For $\lambda\equiv \rho_2=(n-\frac12,n-\frac32,\dots, \frac12)$ we have
\[
{\rm sin}^-_{\rho_2}(x)= \prod_{1\le i<j\le n} \sin \, \pi (x_i-x_j)
\sin \, \pi (x_i+x_j)
 \prod_{1\le i\le n} \sin \pi
x_i.
\]
For $\lambda\equiv \rho_3=(n-1,n-2,\dots, 1,0)$ one has
\[
{\rm sin}^-_{\rho_3}(x)= \prod_{1\le i<j\le n} \sin \, \pi (x_i-x_j)
\sin \, \pi (x_i+x_j).
\]
Similarly, for symmetric multivariate cosine functions we have
\[
{\rm cos}^+_{\rho_1}(x)= \prod_{1\le i<j\le n} \cos \, \pi (x_i-x_j)
\cos \, \pi (x_i+x_j)
 \prod_{1\le i\le n} \cos 2\pi x_i,
\]
\[
{\rm cos}^+_{\rho_2}(x)= \prod_{1\le i<j\le n} \cos \, \pi (x_i-x_j)
\cos \, \pi (x_i+x_j)
 \prod_{1\le i\le n} \cos \pi
x_i,
\]
\[
{\rm cos}^+_{\rho_3}(x)= \prod_{1\le i<j\le n} \cos \, \pi (x_i-x_j)
\cos \, \pi (x_i+x_j).
\]
These formulas are proved in the same way as the formula \eqref{anttttt-3}.

\section{Solutions of the Laplace equation}
The Laplace operator on the $n$-dimensional Euclidean space
$E_n$ in the orthogonal coordinates
$x=(x_1,x_2,\dots,x_n)$ has the
form
\[
\Delta=\frac{\partial^2}{\partial
x^2_1}+\frac{\partial^2}{\partial x^2_2}+\cdots
+\frac{\partial^2}{\partial x^2_n} .
\]

We take any summand in the expression for symmetric or antisymmetric
multivariate sine or cosine function and act upon it by the operator
$\Delta$. We get
\begin{gather*}
\Delta
 \sin 2\pi (w\lambda)_1x_1 \sin 2\pi (w\lambda)_2x_2\cdots  \sin 2\pi
(w\lambda)_n x_{n}
\\
\qquad  {}=-4\pi^2 \langle \lambda,\lambda \rangle\,
\sin 2\pi (w\lambda)_1x_1 \sin 2\pi (w\lambda)_2x_2\cdots  \sin 2\pi
(w\lambda)_n x_{n},
\end{gather*}
where $\lambda=(\lambda_1,\lambda_2,\dots ,\lambda_n)$  determines
${\rm sin}^\pm_\lambda(x)$ or ${\rm cos}^\pm_\lambda(x)$ and
$\langle \lambda,\lambda  \rangle=\sum_{i=1}^n \lambda^2$. Since
this action of $\Delta$ does not depend on a summand from the
expression for symmetric or antisymmetric multivariate sine or
cosine function, we have
\begin{gather}\label{Lapl}
\Delta\, {\rm sin}^\pm_\lambda(x)= -4\pi^2\langle \lambda,\lambda
\rangle\,
 {\rm sin}^\pm_\lambda(x),
 \end{gather}
\begin{gather}\label{Lapl-00}
\Delta\, {\rm cos}^\pm_\lambda(x)= -4\pi^2\langle \lambda,\lambda
\rangle\,
 {\rm cos}^\pm_\lambda(x).
 \end{gather}

The formulas \eqref{Lapl} and \eqref{Lapl-00} admit a generalization.
Let $\sigma_k(y_1,y_2,\dots,y_n)$ be the $k$-th elementary
symmetric polynomial of degree $k$ of the variables
$y_1,y_2,\dots,y_n$, that is,
\[
 \sigma_k(y_1,y_2,\dots,y_n)=\sum_{1\le k_1<k_2<\cdots <k_n\le n}
y_{k_1} y_{k_2}\cdots y_{k_n}.
\]
Then for $k=1,2,\dots,n$ we have
\begin{gather}\label{Lap-gene}
\sigma_k\left( \tfrac{\partial^2}{\partial
x^2_1},\tfrac{\partial^2}{\partial
x^2_2},\dots,\tfrac{\partial^2}{\partial x^2_n} \right) {\rm
sin}^\pm _\lambda(x)= (-4\pi^2)^k \sigma_k
(\lambda_1^2,\lambda_2^2,\dots,\lambda_n^2) {\rm sin}^\pm
_\lambda(x),
 \end{gather}
\[
\sigma_k\left( \tfrac{\partial^2}{\partial
x^2_1},\tfrac{\partial^2}{\partial
x^2_2},\dots,\tfrac{\partial^2}{\partial x^2_n} \right) {\rm
cos}^\pm _\lambda(x)= (-4\pi^2)^k \sigma_k
(\lambda_1^2,\lambda_2^2,\dots,\lambda_n^2) {\rm cos}^\pm
_\lambda(x).
\]
Note that $n$ differential equations \eqref{Lap-gene} are algebraically
independent.

Thus, the functions ${\rm sin}^\pm _m(x)$, ${\rm cos}^\pm _m(x)$,
$m=(m_1m_2,\dots,m_n)$, $m_j\in \mathbb{Z}$, are eigenfunctions of
the operators $\sigma_k\left( \frac{\partial^2}{\partial
x^2_1},\frac{\partial^2}{\partial
x^2_2},\dots,\frac{\partial^2}{\partial x^2_n} \right)$,
$k=1,2,\dots,n$, on the fundamental domain $F(\tilde S_n^{\rm aff})$
satisfying the boundary conditions formulated in section 4. For
example, the functions ${\rm sin}^-_m(x)$ are eigenfunctions of
these operators on the fundamental domain $F(\tilde S_n^{\rm aff})$
satisfying the boundary condition
\begin{gather}\label{Neum1}
{\rm sin}^-_m(x) =0 \qquad {\rm for}\qquad x\in \partial F(\tilde
S^{\rm aff}_n)
\end{gather}
(the Dirichlet boundary value problem). The functions ${\rm
cos}^+_m(x)$ are eigenfunctions of these operators on $F(\tilde
S^{\rm aff}_n)$ satisfying the boudary condition
\[
 \frac{\partial {\rm cos}^+_m(x)}{\partial {\bf n}}=0    \ \ \
{\rm for}\ \ \ x\in \partial F(\tilde S^{\rm aff}_n),
\]
that is, these functions give solutions of the Neumann boudary value
problem on $\partial F(\tilde S^{\rm aff}_n)$.

\section{Symmetric and antisymmetric multivariate sine and cosine
series}
Symmetric and antisymmetric trigonometric functions determine
symmetric and antisymmetric multivariate trigonometric Fourier transforms
which generalize the usual trigonometric Fourier transforms.

As in the case of trigonometric functions of one variable, (anti)symmetric
sine and cosine functions determine three types of
trigonometric Fourier transforms:
 \medskip

(a) Fourier transforms related to the sine and cosine functions
${\rm sin}^\pm_m(x)$ and ${\rm cos}^\pm_m(x)$ with
$m=(m_1,m_2,\dots,m_n)$, $m_j\in \mathbb{Z}$ (trigonometric Fourier
series);
 \medskip

(b) Fourier transforms related to the sine and cosine functions
${\rm sin}^\pm_\lambda(x)$ and ${\rm cos}^\pm_\lambda(x)$ with
$\lambda\in D_+\equiv \overline{D^+_+}$ (integral Fourier
transforms);
 \medskip

(c) Symmetric and antisymmetric multivariate finite sine and
cosine Fourier transforms.
 \medskip

In this section we consider the case (a). Let $f(x)$ be an
antisymmetric (with respect to the extended affine symmetric group
$\tilde S_n^{\rm aff}$) continuous real function on the
$n$-dimensional Euclidean space $E_n$, which has continuous
derivatives and vanishes on the boundary $\partial F(\tilde S^{\rm
aff}_n)$ of the fundamental domain $F(\tilde S^{\rm aff}_n)$, that
is, $f(x)$ behaves under action of elements of $\tilde S^{\rm
aff}_n$ as the functions ${\rm sin}^-_m(x)$ do. We may consider this
function on the set ${\sf T}=[0,\frac12 ]^n$ (this set is a closure
of the union of the sets $w F(\tilde S^{\rm aff}_n)$, $w\in S_n$).
Then $f(x)$, as a function on ${\sf T}$, can be expanded in sine
functions
\[
\sin 2\pi m_1 x_1{\cdot}
\sin 2\pi m_2 x_2\cdots \sin 2\pi m_n x_n,\ \ \
m_i\in \mathbb{Z}^{>0}.
\]
We have
 \begin{equation}\label{expan-anti}
f(x)=\sum_{m_i\in \mathbb{Z}^{>0}} c_m\,
\sin 2\pi m_1 x_1{\cdot}
\sin 2\pi m_2 x_2\cdots \sin 2\pi m_n x_n,
 \end{equation}
where $m=(m_1,m_2,\dots,m_n)$. Let us show that
$c_{wm}=(\det w)c_m$, $w\in S_n$. We represent each sine function
in the expression \eqref{expan-anti} in the form
$\sin \alpha=(2{\rm i})^{-1}(e^{{\rm i}\alpha}-e^{-{\rm i}\alpha})$.
Then
\[
 f(x)=\sum_{m_i\in \mathbb{Z}} c_m e^{2\pi {\rm i}m_1x_1}e^{2\pi {\rm i}m_2x_2}
\cdots e^{2\pi {\rm i}m_nx_n}=\sum_{m_i\in \mathbb{Z}} c_m
e^{2\pi {\rm i}  \langle m, x\rangle},
\]
where $\langle m, x\rangle=\sum_{i=1}^n m_ix_i$ and $c_m$ with positive $m_i$,
$i=1,2,\dots,n$, are
such as in \eqref{expan-anti} and each change of a sign in $m$ leads to
multiplication of $c_m$ by $(-1)$.
Due to the property $f(wx)=(\det w)f(x)$, $w\in S_n$, for any $w\in S_n$ we have
 \begin{align}\label{ddddet}
f(wx)= &
\sum_{m_i\in \mathbb{Z}} c_m
e^{2\pi{\rm i}m_1 x_{w(1)}}
\cdots e^{2\pi{\rm i}m_n x_{w(n)}}
=\sum_{m_i\in \mathbb{Z}} c_{m}
e^{2\pi{\rm i}m_{w^{-1}(1)} x_{1}}
\cdots
e^{2\pi{\rm i}m_{w^{-1}(n)} x_n}
\notag \\
= & \sum_{m_i\in \mathbb{Z}} c_{wm}
e^{2\pi{\rm i}m_1 x_{1}}
\cdots e^{2\pi{\rm i}m_n x_{n}}
=(\det w)f(x)
\notag \\
= & \; (\det w)
\sum_{m_i\in \mathbb{Z}} c_m
e^{2\pi{\rm i}m_1 x_1}
\cdots e^{2\pi{\rm i}m_n x_n}. \notag
\end{align}
Thus, the coefficients $c_m$ in \eqref{expan-anti} satisfy the conditions
$c_{wm}=(\det w)c_m$, $w\in S_n$.

Collecting in \eqref{expan-anti} products of sine functions at
$(\det w) c_{wm}$, $w\in S_n$, we obtain the expansion
 \begin{equation}\label{expan-anti-2}
f(x)= \sum_{ m\in P^+_+} c_{m} \det \left( \sin 2\pi
m_ix_j\right)_{i,j=1}^n \equiv \sum_{m\in P^+_+} c_m\, {\rm
sin}^-_m(x),
 \end{equation}
where $P^+_+:= D^+_+\cap \mathbb{Z}^n$. Thus, {\it any antisymmetric
(with respect to $S_n$) continuous real function $f$ on ${\sf T}$,
which has continuous derivatives, can be expanded in antisymmetric
multivariate sine functions} ${\rm sin}^-_m(x)$, $m\in P^+_+$.

By the orthogonality relation \eqref{symmmm-cos}, the coefficients $c_m$ in the
expansion \eqref{expan-anti-2} are determined by the formula
 \begin{gather}\label{decom-ant-2a}
c_m = 2^{2n} \int_{\overline{F(S^{\rm aff}_n)}} f(x) \det \left(
\sin 2\pi m_ix_j\right)_{i,j=1}^n dx =2^{2n}
\int_{\overline{F(S^{\rm aff}_n)}} f(x) {\rm sin}^-_m(x) dx ,
 \end{gather}
Moreover, the Plancherel formula
\[
\sum_{m\in P_+^+}  |c_m|^2= 2^{2n}
\int_{\overline{F(S^{\rm aff}_n)}} |f(x)|^2dx
 \]
holds, which means that the Hilbert spaces with the appropriate
scalar products are isometric.

Formula \eqref{decom-ant-2a} is an antisymmetrized sine Fourier
transform of the function $f(x)$. Formula \eqref{expan-anti-2} gives
an inverse transform. Formulas \eqref{expan-anti-2} and
\eqref{decom-ant-2a} give the {\it antisymmetric multivariate sine
Fourier transforms} corresponding to antisymmetric sine functions
${\rm sin}^-_m(x)$, $m\in P_+^+$.

Analogous transforms hold for symmetric cosine functions ${\rm
cos}^+_m(x)$, $m\in P_+=D_+\cap \mathbb{Z}^n$. Let $f(x)$ be a
symmetric (with respect to the group $\tilde S_n^{\rm aff}$)
continuous real function on the $n$-dimensional Euclidean space
$E_n$, which has continuous derivatives, that is, $f(x)$ behaves
under action of elements of $\tilde S_n^{\rm aff}$ as the functions
${\rm cos}^+_m(x)$ do. We may consider this function as a function
on $F(\tilde S^{\rm aff}_n)$. Then we can expand this function as
\begin{equation}\label{decom-anti-c}
f(x)=\sum_{{m}\in P_+} c_{m} {\det}^+(\cos 2\pi
m_ix_j)_{i,j=1}^n=\sum_{{m}\in P_+} c_{m} {\rm cos}^+_m(x),
 \end{equation}
where $m=(m_1,m_2,\dots,m_n)$ are {\it integer} $n$-tuples
such that $m_1\ge m_2\ge \cdots \ge m_n\ge 0$. The coefficients
$c_{m}$ of this expansion are given by the formula
 \begin{equation}\label{decom-anti-c-}
 c_{m}=2^{2n} |G_m|^{-1} \int_{\overline{F(S^{\rm aff}_n)}} f(x)
{\rm cos}^+_m(x) dx.
 \end{equation}
The Plancherel formula is of the form
\[
 \sum_{{m}\in P_+} |G_m||c_{ m}|^2=2^{2n} \int_{\overline{F(S^{\rm aff}_n)}}
|f(x)|^2 dx.
\]

Now let  $f(x)$ be an antisymmetric (with respect to the permutation
group $S_n$) continuous real function on the set ${\sf T}=[ 0,\frac12 ]^n$,
which has continuous derivatives and vanishes on the boundary
$\partial F(\tilde S^{\rm aff}_n)$ of the fundamental domain
$F(\tilde S^{\rm aff}_n)$. Then $f(wx)=(\det w)f(x)$, $w\in S_n$.
We consider this function as a
function on $F(\tilde S^{\rm aff}_n)$. One has the expansion
\begin{align}\label{expan-anti-sym-2}
f(x)= \sum_{m\in P_+^+} c_m \det \left( \cos 2\pi
m_ix_j\right)_{i,j=1}^n \equiv \sum_{m\in P_+^+} c_m {\rm
cos}^-_m(x),
 \end{align}
where
 \begin{gather}\label{decom-ant-sym-2}
c_m =2^{2n} \int_{\overline{F(S^{\rm aff}_n)}} f(x) {\rm cos}^-_m(x)
dx.
 \end{gather}
Moreover, the Plancherel formula
$\sum_{m\in P^+_+}  |c_m|^2=2^{2n}
\int_{\overline{F(S^{\rm aff}_n)}} |f(x)|^2dx$ holds.

A similar expansion for the functions ${\rm sin}^+_m(x)$, $m\in
P_+$, is of the form
\begin{equation}\label{decom-anti-cc}
f(x)=\sum_{{m}\in P_+} c_{m} {\det}^+(\sin 2\pi
m_ix_j)_{i,j=1}^n\equiv \sum_{{m}\in P_+} c_{m} {\rm sin}^+_m(x),
 \end{equation}
where the coefficients $c_{m}$ are given by
 \begin{equation}\label{decom-anti-cc-}
 c_{m}=2^{2n} |G_m|^{-1}\int_{\overline{F(S^{\rm aff}_n)}} f(x)
 {\rm sin}^+_m(x) dx.
 \end{equation}
The Plancherel formula is of the form
$\sum_{{m}\in P_+} |G_m| |c_{m}|^2=2^{2n} \int_{\overline{F(S^{\rm aff}_n)}}
|f(x)|^2 dx$.

\section{Symmetric and antisymmetric multivariate sine and cosine\\
Fourier transforms on $F(\tilde S_n)$}
The expansions of the previous subsection give expansions of
functions on the fundamental domain $F(\tilde S^{\rm aff}_n)$ in
functions ${\rm sin}^\pm_m(x)$ and ${\rm cos}^\pm_m(x)$ with
integral $m=(m_1,m_2,\dots,m_n)$.  The functions ${\rm
sin}^\pm_\lambda(x)$ and ${\rm cos}^\pm_\lambda(x)$ with $\lambda$
lying in the fundamental domain $F(\tilde S_n)$ (and not obligatory
integral) are invariant (anti-invariant) only with respect to the
extended permutation group $\tilde S_n$. A fundamental domain of
$\tilde S_n$ coincides with the set $D^+_+$ from section 3. The sine
and cosine functions ${\rm sin}^\pm_\lambda(x)$ and ${\rm
cos}^\pm_\lambda(x)$, determined by $\lambda\in D_+$, give Fourier
transforms on the domain $D_+$.

We began with the usual sine Fourier transforms on ${\mathbb R}_+^n$:
 \begin{gather}\label{F-1}
\tilde f (\lambda)=\int_{{\mathbb R}^n_+} f(x) \sin 2\pi \lambda_1 x_1\,
\sin 2\pi \lambda_2 x_2 \cdots \sin 2\pi \lambda_n x_n\, dx,
  \\
  \label{F-2}
 f (x)=2^{2n}\int_{{\mathbb R}_+^n} \tilde f(\lambda)
\sin 2\pi \lambda_1 x_1\,
\sin 2\pi \lambda_2 x_2 \cdots \sin 2\pi \lambda_n x_n\, d\lambda.
  \end{gather}
Let the function $f(x)$, given on $\mathbb{R}^n_+$, be
anti-invariant with respect to the symmetric group $S_n$, that is,
$f(wx)=(\det w)f(x)$, $w\in S_n$. It is easy to check that the
function $\tilde f (\lambda)$ is also anti-invariant with respect to
the group $S_n$. Replace in \eqref{F-1} $\lambda$ by $w\lambda$,
$w\in S_n$, multiply both sides by $\det w$, and sum these both side
over $w\in S_n$. Due to the expression \eqref{det-B} for symmetric
sine functions ${\rm sin}^-_\lambda(x)$, instead of \eqref{F-1} we
obtain
 \begin{gather}\label{F-3}
\tilde f (\lambda)= \vert S_n\vert^{-1} \int_{{\mathbb R}_+^n} f(x)
{\rm sin}^-_\lambda(x) dx \equiv \int_{D_+} f(x) {\rm
sin}^-_\lambda(x) dx,\qquad \lambda\in D^+_+,
  \end{gather}
where we have taken into account that $f(x)$ is anti-invariant
with respect to $S_n$.

Starting from \eqref{F-2}, we obtain the inverse
formula,
 \begin{gather}\label{F-4}
 f (x)= 2^{2n}  \int_{D_+} \tilde f(\lambda)
{\rm sin}^-_\lambda(x) d\lambda .
  \end{gather}
For the transforms \eqref{F-3} and \eqref{F-4} the Plancherel
formula
 \[
 \int_{D_+} |f(x)|^2 dx=
2^{2n} \int_{D_+} |\tilde f(\lambda) |^2  d\lambda
 \]
holds. The formulas \eqref{F-3} and \eqref{F-4} determine the {\it
antisymmetric multivariate sine Fourier transforms on the domain}
$F(\tilde S_n)$.

Similarly, starting from formulas
\eqref{F-1} and \eqref{F-2} we receive the symmetric multivariate
sine Fourier transforms on the domain $F(\tilde S_n)$:
 \begin{gather}\label{F-5}
\tilde f (\lambda)= \int_{D_+} f(x) {\rm sin}^+_\lambda(x) dx,\qquad
\lambda\in D_+,
  \end{gather}
 \begin{gather}\label{F-6}
 f (x)=2^{2n}
 \int_{D_+} \tilde f(\lambda)
 {\rm sin}^+_\lambda(x) d\lambda .
  \end{gather}
The corresponding Plancherel formula holds.

The cosine functions $ {\rm cos}^\pm_\lambda(x)$ determine similar
transforms. Namely, we have
 \begin{equation}\label{anty-1}
\tilde f (\lambda)= \int_{D_+} f(x) \det \left( \cos
2\pi\lambda_ix_j\right)_{i,j=1}^n  dx \equiv   \int_{D_+} f(x) {\rm
cos}^-_\lambda(x) dx,
  \end{equation}
where
 \begin{equation}\label{anty-2}
 f (x)=2^{2n} \int_{D_+} \tilde f(\lambda)
\det \left( \cos 2\pi\lambda_ix_j\right)_{i,j=1}^n  d\lambda \equiv
2^{2n} \int_{D_+} \tilde f(\lambda) {\rm cos}^-_\lambda(x) d\lambda,
  \end{equation}
and
 \begin{equation}\label{anty-3}
\tilde f (\lambda)= \int_{D_+} f(x) {\det}^+ \left( \cos
2\pi\lambda_ix_j\right)_{i,j=1}^n  dx \equiv \int_{D_+} f(x) {\rm
cos}^+_\lambda(x) dx,
  \end{equation}
where
 \begin{equation}\label{anty-4}
 f (x)=2^{2n} \int_{D_+} \tilde f(\lambda)
{\det}^+ \left( \cos 2\pi\lambda_ix_j\right)_{i,j=1}^n  d\lambda
\equiv 2^{2n} \int_{D_+} \tilde f(\lambda) {\rm cos}^+_\lambda(x)
d\lambda .
  \end{equation}
The corresponding Plancherel formulas hold.

\section{Finite 1-dimensional sine and cosine transforms}
Finite one-dimensional sine and cosine transforms are useful for
applications. The theory of these transforms as well as
their different applications and methods of work with them
are given in Ref. 17. In this
section we give these one-dimensional transforms in the
form$^{11}$ which will be used in the following sections.

Let $N$ be a positive integer. To this number there corresponds the
finite set of points (the grid) $\frac rN$, $r=0,1,2,\dots,N$. We
denote this set by $F_N$,
 \begin{equation}\label{grid-3}    \textstyle{
F_N=\left\{ 0,\frac1N, \frac2N,\dots,\frac{N-1}N, 1\right\} .}
 \end{equation}

We consider sine functions on the grid $F_N$, that is, the functions
 \begin{equation}\label{grid-4}
\varphi_m(s):=\sin (\pi ms) ,\ \ \ s\in F_N,\ \ m\in {\mathbb
Z}^{\ge}.
 \end{equation}

Since $\varphi_m(s)=\pm \varphi_{m+N}(s)$ and $\varphi_0(s)=\varphi_N(s)=0$,
we consider these discrete functions only for
\[
m\in D_N:=\{ 1,2,\dots,N-1\}.
\]
The functions \eqref{grid-4} vanish on the points 0 and 1 of $F_N$. For this reason,
we also consider the subset
\[      \textstyle{
 F^-_N=\left\{ \frac1N, \frac2N,\dots,\frac{N-1}N\right\}
 \ \ \ (N{-}1\ {\rm points}) }
\]
of the grid $F_N$.

The functions \eqref{grid-4} are orthogonal on $F^-_N$ and
the orthogonality relation is of the form
 \begin{equation}\label{grid-44}
 \langle \varphi_m,\varphi_{m'}  \rangle=\sum_{s\in F_N^-}
\varphi_m(s) \varphi_{m'}(s)=\frac N2 \delta_{mm'},\qquad m,m'\in
D_N.
\end{equation}
They determine the following expansion of
functions, given on the grid $F^-_N$:
 \begin{equation}\label{grid-5}
f(s)=\sum_{m=1}^{N-1} a_m \sin (\pi ms) ,
 \end{equation}
where the coefficients $a_m$ are given by
 \begin{equation}\label{grid-6}
a_m=\frac2{N} \sum_{s\in F_N^-} f(s)\sin (\pi ms).
 \end{equation}
Formulas \eqref{grid-5} and \eqref{grid-6} determine the {\it
discrete sine transform}.

We also consider cosine functions on the grid $F_N$, that is, the functions
 \begin{equation}\label{grid-12}
 \phi_m(s)=\cos (\pi ms),\ \ \ s\in F_N,\ \ m\in \{ 0,1,2,\dots N\} .
 \end{equation}
These functions are orthogonal on the grid $F_N$ with the
orthogonality relation
 \begin{equation}\label{grid-13}
 \langle \phi_m,\phi_{m'}  \rangle=\sum_{s\in F_N} c_s
\phi_m(s)\phi_{m'}(s)=r_m N\delta_{mm'},
 \end{equation}
where $r_m=1$ for $m=0,N$ and $r_m=\frac12$ otherwise, $c_s=\frac 12$ for
$s=0,\, 1$ and $c_s=1$ otherwise${}^{18}$.

The functions \eqref{grid-12} determine an
expansion of functions on the grid $F_N$ as
 \begin{equation}\label{grid-7}
f(s)=\sum_{m=0}^{N} b_m \cos (\pi ms) ,\ \ \ s\in F_N,
 \end{equation}
where the coefficients $b_m$ are given by
 \begin{equation}\label{grid-8}
b_m=r_m^{-1}N^{-1} \sum_{s\in F_N} c_s f(s) \cos (\pi ms).
 \end{equation}
Formulas \eqref{grid-7} and \eqref{grid-8} determine the {\it
discrete cosine transform}.

\section{Antisymmetric multivariate finite sine transforms}
The finite sine and cosine transforms of the previous section can be
generalized to the $n$-dimensional case in symmetric and
antisymmetric forms. In fact, these generalizations are finite
(anti)symmetric multivariate trigonometric transforms. They are
derived by using 1-dimensional finite sine and cosine transforms.
Some of the transforms can be also derived by using the results of
Ref. 18. In order to introduce multivariate finite sine transforms
we have to define (anti)symmetric multivariate finite sine
functions. Note that notations ${\rm sin}^\pm_{\bf m}({\bf s})$ in
this section slightly differ from notations of section 2.

We take the discrete sine function \eqref{grid-4} and make a
multivariate discrete sine function by multiplying $n$ copies of
this function:
\begin{equation}\label{multi-1}
{\rm sin}_{\bf m}({\bf s}):= \sin (\pi m_1s_1) \sin (\pi
m_2s_2)\cdots \sin (\pi m_ns_n),
 \end{equation}
\[
 s_j\in F_N,\ \ \ \ m_i\in D_N\equiv \{ 1,2,\dots,N-1\},
\]
where ${\bf s}=(s_1,s_2,\dots,s_n)$ and ${\bf
m}=(m_1,m_2,\dots,m_n)$. We take these multivariate functions
for integers $m_i$ such that $N>m_1>m_2>\cdots>m_n>0$ and make an
antisymmetrization to obtain a finite version of the
antisymmetric multivariate sine function \eqref{det-B}:
\begin{equation}\label{multi-2}
{\rm sin}^-_{\bf m}({\bf s}):= |S_n|^{-1/2} \det (\sin \pi
m_is_j)_{i,j=1}^n,
 \end{equation}
where $|S_n|$ is an order of the symmetric group $S_n$. (We have
here expressions $\sin \pi m_is_j$, not $\sin 2\pi m_is_j$ as in
\eqref{det-B}.) Since functions $\sin \pi m_is_j$ are considered
for positive $m_i$ and $s_j$, we deal here with the permutation
group $S_n$ instead of the group $\tilde S_n$.

The $n$-tuple ${\bf s}$ in \eqref{multi-2} runs over $( F^-_N)^n
\equiv F^-_N\times \cdots \times F^-_N$ ($n$ times). We
denote by $\hat F_N^n$ the subset of $(F^-_N)^n$ consisting of
${\bf s}\in (F^-_N)^n$ such that
\[
 s_1>s_2>\cdots >s_n.
\]
Note that $s_i$ may take the values $\frac1N,\frac2N,\dots,
\frac{N-1}N$. Acting by permutations $w\in S_n$ upon $\hat F_N^n$ we
obtain the whole set $(F^-_N)^n$ without those points which are
invariant under some nontrivial permutation $w\in S_n$.
Due to antisymmetricity, the
functions \eqref{multi-2} vanishes on the last points.

We denote by $\hat D_N^n$ the set of integer $n$-tuples ${\bf
m}=(m_1,m_2,\dots,m_n)$ such that
\[
 N>m_1>m_2>\cdots>m_n>0.
\]

We need a scalar product of functions \eqref{multi-2}. For
this we define a scalar product of functions \eqref{multi-1} as
\[
\langle \sin_{\bf m}({\bf s}), \sin_{{\bf m}'}({\bf s})
\rangle =\prod_{i=1}^n \langle \varphi_{m_i}(s_i),
\varphi_{m'_i}(s_i)\rangle=\left( \frac N2 \right) ^n \delta_{{\bf m},{\bf m}'},
\]
where the scalar product $\langle \varphi_{m_i}(s_i),
\varphi_{m'_i}(s_i)\rangle$ is given by formula \eqref{grid-44}.
Since functions ${\rm sin}^-_{\bf m}({\bf s})$ are linear
combinations of functions $\sin_{{\bf m}'}({\bf s})$, a scalar
product for ${\rm sin}^-_{\bf m}({\bf s})$ is also defined.

\medskip

\noindent {\bf Proposition 1.} {\it For ${\bf m},{\bf m}'\in
\hat D_N^n$, the discrete functions \eqref{multi-2} satisfy the
orthogonality relation}
\[
\langle {\rm sin}^-_{\bf m}({\bf s}), {\rm sin}^-_{{\bf m}'}({\bf
s})\rangle = \sum_{{\bf s}\in (F^-_N)^n} {\rm sin}^-_{\bf m}({\bf
s}) {\rm sin}^-_{{\bf m}'}({\bf s}) = |S_n| \sum_{{\bf s}\in \hat
F_N^n} {\rm sin}^-_{\bf m}({\bf s}) {\rm sin}^-_{{\bf m}'}({\bf s})
\]
\begin{equation}\label{multi-3}
=(N/2)^n \delta_{{\bf m}{\bf m}'}.
 \end{equation}

\noindent {\bf Proof.} Since $N>m_1>m_2>\cdots>m_n>0$
for ${\bf m}\in \hat D^n_N$, then
due to the orthogonality relation \eqref{grid-44}
for the sine functions $\sin (\pi ms)$ we have
\[
\sum_{{\bf s}\in (F^-_N)^n} {\rm sin}^-_{\bf m}({\bf s}) {\rm
sin}^-_{{\bf m}'}({\bf s})
 = |S_n|^{-1} \sum_{w\in
S_n} \prod_{i=1}^n \sum_{s_i=1}^{N-1} \sin (\pi m_{w(i)}s_i) \sin
(\pi m'_{w(i)}s_i)
\] \[
 = (N/2)^n \delta_{{\bf m}{\bf m}'},
\]
where $(m_{w(1)},m_{w(2)},\dots ,m_{w(n)})$ is obtained from
$(m_1,m_2,\dots, m_n)$ by action by the permutation $w\in S_n$.
Since functions ${\rm sin}^-_{\bf m}({\bf s})$ are antisymmetric
with respect to $S_n$, we have
\[
 \sum_{{\bf s}\in (F^-_N)^n}
{\rm sin}^-_{\bf m}({\bf s}) {\rm sin}^-_{{\bf m}'}({\bf s})
 =|S_n| \sum_{{\bf s}\in \hat
F_N^n} {\rm sin}^-_{\bf m}({\bf s}) {\rm sin}^-_{{\bf m}'}({\bf s}).
\]
This proves the proposition.
 \medskip

Let $f$ be a function on $\hat F^n_N$ (or an antisymmetric function
on $(F^-_N)^n$). Then it can be expanded in functions
\eqref{multi-2} as
\begin{equation}\label{multi-4}
f({\bf s})=\sum_{{\bf m}\in \hat D_N^n}a_{\bf m} {\rm sin}^-_{\bf
m}({\bf s}),
 \end{equation}
where the coefficients $a_{\bf m}$ are determined by the formula
\begin{equation}\label{multi-5}
a_{\bf m}=(2/N)^{n}|S_n| \sum_{{\bf s}\in \hat F_N^n} f({\bf s})
{\rm sin}^-_{\bf m}({\bf s}).
 \end{equation}

A validity of the expansions \eqref{multi-4} and \eqref{multi-5}
follows from the facts that numbers of elements in $\hat D_N^n$ and in
$\hat F_N^n$ are the same and from the orthogonality relation
\eqref{multi-3}.

\section{Symmetric multivariate finite cosine transforms}
We take the finite cosine functions \eqref{grid-12} and make
multivariate finite cosine functions by multiplying $n$ copies of
this function:
\begin{equation}\label{multiv-1}
\cos_{\bf m}({\bf s}):=\cos (\pi m_1s_1) \cos (\pi m_2s_2)\cdots \cos
(\pi m_ns_n),
 \end{equation}
\[
 s_j\in F_N ,\ \ \ \ m_i\in \{ 0,1,2,\dots,N\} .
\]
We consider these functions for integers $m_i$ such that $N\ge m_1\ge
m_2\ge \cdots\ge m_n\ge 0$ (the collection of these $n$-tuples ${\bf
m}=(m_1,m_2,\dots,m_n)$ will be denoted by $\breve D_N^n$)
and make a symmetrization. As a result,
we obtain a finite version of the symmetric multivariate
cosine function \eqref{det-B-C}:
\begin{equation}\label{multiv-2}
{\rm cos}^+_{\bf m}({\bf s}):= |S_n|^{-1/2} \sum_{w\in S_n} \cos
(\pi m_{w(1)}s_1) \cos (\pi m_{w(2)}s_2)\cdots \cos (\pi
m_{w(n)}s_n).
 \end{equation}
(We have here expressions $\cos \pi m_is_j$, not $\cos 2\pi m_is_j$
as in \eqref{det-B-C}. Therefore, the notation ${\rm cos}^+_{\bf
m}({\bf s})$ here slightly differs from the notation in section 2.)

The $n$-tuple ${\bf s}$ in \eqref{multiv-2} runs over $F_N^n$.
We denote by $\breve F_N^n$ the subset of $F_N^n$
consisting of ${\bf s}\in F_N^n$ such that
\[
 s_1\ge s_2\ge \cdots \ge s_n .
\]
Note that $s_i$ here may take the values $0,\frac1N,\frac2N,\dots,
\frac{N-1}N,1$. Acting by permutations $w\in S_n$ upon $\breve
F_N^n$ we obtain the whole set $F_N^n$, where points, invariant
under some nontrivial permutation $w\in S_n$, are repeated several
times. It is easy to see that a point ${\bf s}_0\in F_N^n$ is
repeated $|S_{{\bf s}_0}|$ times in the set $\{ w\breve F_N^n$;
$w\in S_n\}$, where $|S_{{\bf s}_0}|$ is an order of the subgroup
$S_{{\bf s}_0}\subset S_n$, whose elements leave ${\bf s}_0$
invariant.

A scalar product of functions \eqref{multiv-1} is determined by
\[
\langle \cos_{\bf m}({\bf s}),\cos_{{\bf m}'}({\bf s})
\rangle =\prod_{i=1}^n \langle \cos_{m_i}(s_i),
\cos_{m'_i}(s_i)\rangle=N^n r_{m_1}\cdots r_{m_n} \delta_{{\bf m},{\bf m}'},
\]
where the scalar product $\langle \cos_{m_i}(s_i),
\cos_{m'_i}(s_i)\rangle$ is given by \eqref{grid-13}. Since
functions ${\rm cos}^+_{\bf m}({\bf s})$ are linear combinations of
functions $\cos_{{\bf m}'}({\bf s})$, then a scalar product for
${\rm cos}^+_{\bf m}({\bf s})$ is also defined.
 \medskip

\noindent {\bf Proposition 2.} {\it For ${\bf m},{\bf m}'\in \breve
D_N^n$, the discrete functions \eqref{multiv-2} satisfy the
orthogonality relation
\begin{alignat}{2}\label{multiv-3}
\langle {\rm cos}^+_{\bf m}({\bf s}), {\rm cos}^+_{{\bf m}'}({\bf
s})\rangle =&\; \sum_{{\bf s}\in  F_N^n} c_{\bf s} {\rm cos}^+_{\bf
m}({\bf s})
{\rm cos}^+_{{\bf m}'}({\bf s}) \notag\\
=&\; |S_n|\sum_{{\bf s}\in \breve F_N^n} |S_{\bf s}|^{-1} c_{\bf s}
{\rm cos}^+_{\bf m}({\bf s})
{\rm cos}^+_{{\bf m}'}({\bf s}) \notag\\
=&\; N^n  r_{\bf m} |S_{\bf m}|  \delta_{{\bf m}{\bf m}'},
\end{alignat}
where $c_{\bf s}=c_{s_1}c_{s_2}\cdots c_{s_n}$, $r_{\bf
s}=r_{m_1}r_{m_2}\cdots r_{m_n}$, and $c_{s_i}$ and  $r_{m_i}$ are
such as in formula \eqref{grid-13}.}
\medskip

\noindent {\bf Proof.} Due to the orthogonality relation for the
cosine functions $\phi_m(s)=\cos (\pi ms)$ (see formula
\eqref{grid-13}) we have
\begin{alignat}{2}\label{multiv-12}
\sum_{{\bf s}\in F_N^n} c_{\bf s} {\rm cos}^+_{\bf m}({\bf s}) {\rm
cos}^+_{{\bf m}'}({\bf s}) =&\;
 \frac{|S_{\bf m}|}{|S_n|^{-1}|} \sum_{w\in S_n} \prod_{i=1}^n
\sum_{s_i=0}^{N}c_{s_i} \cos (\pi m_{w(i)}s_i) \cos (\pi m'_{w(i)}s_i)\notag \\
=&\; |S_{\bf m}| N^n r_{\bf m} \delta_{{\bf m}{\bf m}'},
\end{alignat}
where $(m_{w(1)},m_{w(2)},\dots ,m_{w(n)})$ is obtained from
$(m_1,m_2,\dots, m_n)$ by action by the permutation $w\in S_n$.
Since functions ${\rm cos}^+_{\bf m}({\bf s})$ are symmetric with
respect to $S_n$, we have
\[
 \sum_{{\bf s}\in F_N^n} c_{\bf s}
{\rm cos}^+_{\bf m}({\bf s}) {\rm cos}^+_{{\bf m}'}({\bf s}) =|S_n|
\sum_{{\bf s}\in \breve F_N^n} |S_{\bf s}|^{-1} c_{\bf s} {\rm
cos}^+_{\bf m}({\bf s}) {\rm cos}^+_{{\bf m}'}({\bf s}) .
\]
This proves the proposition.
 \medskip

Let $f$ be a function on $\breve F^n_N$ (or a symmetric
function on $F_N^n$). Then it can be expanded in functions
\eqref{multiv-2} as
\begin{equation}\label{multiv-4}
f({\bf s})=\sum_{{\bf m}\in \breve D_N^n}a_{\bf m} {\rm cos}^+_{\bf
m}({\bf s}),
 \end{equation}
where the coefficients $a_{\bf m}$ are determined by the formula
\begin{alignat}{2}\label{multiv-5}
a_{\bf m}=&\; N^{-n}|S_{\bf m}|^{-1} r_{\bf m}^{-1} \langle f({\bf
s}), {\rm cos}^+_{\bf m}({\bf s})  \rangle
\notag\\
=&\; N^{-n} |S_{\bf m}|^{-1} r_{\bf m}^{-1} |S_n|\sum_{{\bf s}\in
\breve F_N^n} |S_{\bf s}|^{-1} c_{\bf s}f({\bf s}) {\rm cos}^+_{\bf
m}({\bf s}).
 \end{alignat}

A validity of the expansions \eqref{multiv-4} and \eqref{multiv-5}
follows from the fact that numbers of elements in $\breve  D_N^n$
and $\breve  F_N^n$ are the same and from the orthogonality relation
\eqref{multiv-3}.

\section{Other 1-dimensional finite cosine transforms}
Along with the finite cosine transform of section 9 there
exist other 1-dimensional finite transforms with the discrete cosine function as
a kernel (see, for example, Refs. 19 and 20).
In Ref. 19 the finite cosine transforms are denoted
as DCT-1, DCT-2, DCT-3, DCT-4. The transform DCT-1 is in fact the
transform, considered in section 9. Let us expose all these
transforms (including the transform DCT-1), conserving
notations used in the literature on signal processing. They are determined by a
positive integer $N$.
\medskip

\noindent {\bf DCT-1.} This transform is given by the kernel
\[
\mu_r(k)=\cos \frac{\pi rk}N,  \qquad {\rm where}\qquad
k,r\in \{ 0,1,2,\dots, N\}.
\]
The orthogonality relation for these discrete functions is given by
 \begin{equation}\label{DCT-1-3}
\sum_{k=0}^N c_k  \cos \frac{\pi rk}N  \cos \frac{\pi r'k}N =h_r
\frac N2 \delta_{rr'},
 \end{equation}
where $c_k=\frac 12$ for $k=0,N$ and $c_k=1$ otherwise, $h_r=2$ for
$r=0,N$ and $h_r=1$ otherwise.

Thus, these functions give the expansion
 \begin{equation}\label{DCT-1-4}
f(k)=\sum_{r=0}^N a_r  \cos \frac{\pi rk}N ,\qquad {\rm where}\qquad
 a_r= \frac2{h_rN} \sum_{k=0}^N c_kf(k) \cos \frac{\pi rk}N .
 \end{equation}

\noindent {\bf DCT-2.} This transform is given by the kernel
 \[
\omega_r(k)=\cos \frac{\pi (r+\frac12)k}N, \qquad {\rm where}\qquad
k,r\in \{ 0,1,2,\dots, N-1\}.
\]
The orthogonality relation for these discrete functions is given by
 \begin{equation}\label{DCT-2-3}
\sum_{k=0}^{N-1} c_k  \cos \frac{\pi (r+\frac12)k}N \cos \frac{\pi
(r'+\frac12)k}N= \frac N2 \delta_{rr'},
 \end{equation}
where $c_k=1/2$ for $k=0$ and $c_k=1$ otherwise.

These functions determine the expansion
 \begin{equation}\label{DCT-2-4}
f(k)=\sum_{r=0}^{N-1} a_r  \omega_r(k), \qquad
 {\rm where} \qquad
 a_r= \frac2{N} \sum_{k=0}^{N-1} c_kf(k) \omega_r(k).
 \end{equation}

\noindent {\bf DCT-3.} This transform is determined by the kernel
 \[
\sigma_r(k)=\cos \frac{\pi r(k+\frac12)}N,
\]
where $k$ and $r$ run over the values $\{ 0,1,2,\dots,N-1\}$. The
orthogonality relation for these discrete functions is given by the
formula
 \begin{equation}\label{DCT-3-3}
\sum_{k=0}^{N-1}  \cos \frac{\pi r(k+\frac12)}N \cos \frac{\pi
r'(k+\frac12)}N=h_r \frac N2 \delta_{rr'},
 \end{equation}
where $h_k=2$ for $k=0$ and $h_k=1$ otherwise.

These functions give the expansion
 \begin{equation}\label{DCT-3-4}
f(k)=\sum_{r=0}^{N-1} a_r  \cos \frac{\pi r(k+\frac12)}N, \quad
 {\rm where} \quad
 a_r= \frac2{h_rN} \sum_{k=0}^{N-1} f(k) \cos \frac{\pi r(k+\frac12)}N.
 \end{equation}

\noindent {\bf DCT-4.} This transform is given by the kernel
 \[
\tau_r(k)=\cos \frac{\pi (r+\frac12)(k+\frac12)}N,
 \]
where $k$ and $r$ run over the values $\{ 0,1,2,\dots,N-1\}$. The
orthogonality relation for these discrete functions is given by
 \begin{equation}\label{DCT-4-3}
\sum_{k=0}^{N-1}  \cos \frac{\pi (r+\frac12)(k+\frac12)}N \cos
\frac{\pi (r'+\frac12)(k+\frac12)}N= \frac N2 \delta_{rr'}.
 \end{equation}
These functions determine the expansion
 \begin{equation}\label{DCT-4-4}
f(k)=\sum_{r=0}^{N-1} a_r  \cos \frac{\pi (r+\frac12)(k+\frac12)}N,
\;\;
 {\rm where} \;\;
 a_r= \frac2{N} \sum_{k=0}^{N-1} f(k) \cos \frac{\pi (r+\frac12)(k+\frac12)}N.
 \end{equation}

Note that there exist also four discrete sine transforms,
corresponding to the above discrete cosine transforms. They are
obtained from the cosine transforms by replacing in \eqref{DCT-1-4},
\eqref{DCT-2-4}, \eqref{DCT-3-4} and \eqref{DCT-4-4} cosines
discrete functions by sine discrete functions${}^{17, 20}$.

\section{Other antisymmetric multivariate finite cosine transforms}
Each of the finite cosine transforms DCT-1, DCT-2, DCT-3, DCT-4
generates the corresponding antisymmetric multivariate finite
cosine transform. We call them AMDCT-1,  AMDCT-2,
AMDCT-3 and AMDCT-4. Let us give these transforms without proof.
Their proofs are the same as in the case of symmetric
multivariate finite cosine transforms of section 11. Below we
use the notation $\tilde D^{n}_N$ for the subset of the set $D_N^n\equiv
D_N\times D_N\times \dots \times D_N$ ($n$ times) with $D_N=\{
0,1,2,\dots,N\}$ consisting of points ${\bf r}=(r_1,r_2,\dots,r_n)$,
$r_i\in D_N$, such that
\[
 N\ge r_1>r_2>\cdots >r_n\ge 0.
\]

{\bf AMDCT-1.}  This transform is given by the kernel
 \begin{equation}\label{AMDCT-1-1}
{\rm cos}^-_{\bf r}({\bf k})\equiv {\rm cos}^{1,-}_{\bf r}({\bf
k})=|S_n|^{-1/2} \det \left( \cos \frac{\pi
r_ik_j}N\right)_{i,j=1}^n,
 \end{equation}
where ${\bf r}\in \tilde D^{n}_N$ and ${\bf
k}=(k_1,k_2,\dots,k_n)$, $k_i\in \{ 0,1,2,\dots,N\}$. The
orthogonality relation for these kernels is
 \begin{equation}\label{AMDCT-1-2}
\langle {\rm cos}^-_{\bf r}({\bf k}), {\rm cos}^-_{{\bf r}'}({\bf
k}) \rangle = |S_n|\sum_{{\bf k}\in \tilde D^{n}_N}c_{\bf k} {\rm
cos}^-_{\bf r}({\bf k}) {\rm cos}^-_{{\bf r}'}({\bf k}) =h_{\bf r}
\left( \frac N2\right)^n \delta_{{\bf r}{\bf r}'},
 \end{equation}
where
\[
 c_{\bf k}=c_1c_2\cdots c_n,\ \ \  h_{\bf k}=h_1h_2\cdots h_n,
\]
and $c_i$ and $h_j$ are such as in formula \eqref{DCT-1-3}.

This transform is given by the formula
 \begin{equation}\label{AMDCT-1-3}
f( {\bf k})=\sum _{{\bf r}\in \tilde D^{n}_N} a_{\bf r}
 {\rm cos}^-_{\bf r}({\bf k}) , \ \ {\rm where} \ \
a_{\bf r}= h_{\bf r}^{-1}|S_n|\left( \frac2N\right)^n \sum_{{\bf
k}\in \tilde D^{n}_N} c_{\bf k} f({\bf k}) {\rm cos}^-_{\bf r}({\bf
k}).
 \end{equation}
The corresponding Plancherel formula is
\[
 |S_n| \sum_{{\bf k}\in \tilde D^{n}_N} c_{\bf k} |f({\bf k})|^2=
\left( \frac N2\right)^n \sum _{{\bf r}\in \tilde D^{n}_N}  h_{\bf r}
|a_{\bf r}|^2.
\]

{\bf AMDCT-2.}  Let $\tilde D_{N-1}^{n}$ be the subset of
$D^n_{N-1}$ (with $D_{N-1}=\{ 0,1,\dots,N-1\}$) consisting of points
${\bf r}=(r_1,r_2,\dots,r_n)$, $r_i\in D_{N-1}$, such that
\[
 N-1\ge r_1>r_2>\cdots >r_n\ge 0.
\]
This transform is given by the kernel
 \begin{equation}\label{AMDCT-2-1}
{\rm cos}^{2,-}_{\bf r}({\bf k})=|S_n|^{-1/2}\det \left( \cos
\frac{\pi (r_i+\frac12)k_j}N\right)_{i,j=1}^n,
 \end{equation}
where ${\bf r}\in \tilde D^{n}_N$ and ${\bf k}=
(k_1,k_2,\dots,k_n)$, $k_i\in \{ 0,1,2,\dots,N-1\}$. The
orthogonality relation for these kernels is
 \begin{equation}\label{AMDCT-2-2}
\langle {\rm cos}^{2,-}_{\bf r}({\bf k}),{\rm cos}^{2,-}_{{\bf
r}'}({\bf k}) \rangle = |S_n|\sum_{{\bf k}\in \tilde
D_{N-1}^{n}}c_{\bf k} {\rm cos}^{2,-}_{\bf r}({\bf k})\,{\rm
cos}^{2,-}_{{\bf r}'}({\bf k}) = \left(\frac N2\right)^n
\delta_{{\bf r}{\bf r}'},
 \end{equation}
where $c_{\bf k}=c_1c_2\cdots c_n$ and $c_i$ are such as in formula
\eqref{DCT-2-3}.

This transform is given by the formula
 \begin{equation}\label{AMDCT-2-3}
f( {\bf k}){=}\sum _{{\bf r}\in \tilde D_{N-1}^{n}} a_{\bf r} {\rm
cos}^{2,-}_{\bf r}({\bf k}) , \quad {\rm where} \quad a_{\bf
r}{=}|S_n| \left( \frac2N\right)^n \sum_{{\bf k}\in \tilde
D_{N-1}^{n}} c_{\bf k} f({\bf k}) {\rm cos}^{2,-}_{\bf r}({\bf k}).
 \end{equation}
The Plancherel formula is of the form
\[
 |S_n| \sum_{{\bf k}\in \tilde D_{N-1}^{n}} c_{\bf k} |f({\bf k})|^2=
\left( \frac N2\right)^n \sum _{{\bf r}\in \tilde D_{N-1}^{n}} |a_{\bf
r}|^2.
\]

{\bf AMDCT-3.}  This transform is given by the kernel
 \begin{equation}\label{AMDCT-3-1}
{\rm cos}^{3,-}_{\bf r}({\bf k})=|S_n|^{-1/2} \det \left( \cos
\frac{\pi r_i(k_j+\frac12)}N\right)_{i,j=1}^n,\ \ \ {\bf r}\in
\tilde D_{N-1}^{n},\ \ \ k_j\in D_{N-1}.
 \end{equation}
The orthogonality relation for these kernels is
 \begin{equation}\label{AMDCT-3-2}
\langle {\rm cos}^{3,-}_{\bf r}({\bf k}),{\rm cos}^{3,-}_{{\bf
r}'}({\bf k}) \rangle = |S_n|\sum_{{\bf k}\in \tilde D_{N-1}^{n}}
{\rm cos}^{3,-}_{\bf r}({\bf k}) {\rm cos}^{3,-}_{{\bf r}'}({\bf k})
=h_{\bf r} \left(\frac N2\right)^n \delta_{{\bf r}{\bf r}'},
 \end{equation}
where $h_{\bf r}=h_1h_2\cdots h_n$ and $h_j$ are such as in formula
\eqref{DCT-3-3}.

This transform is given by the formula
 \begin{equation}\label{AMDCT-3-3}
f( {\bf k})=\sum _{{\bf r}\in \tilde D_{N-1}^{n}} a_{\bf r}
 {\rm cos}^{3,-}_{\bf r}({\bf k}) ,\quad {\rm where} \quad
a_{\bf r}=\frac{|S_n|}{h_{\bf r}} \left( \frac2N\right)^n \sum_{{\bf
k}\in \tilde D_{N-1}^{n}}  f({\bf k}) {\rm cos}^{3,-}_{\bf r}({\bf
k}).
 \end{equation}
The Plancherel formula is of the form
\[
 |S_n| \sum_{{\bf k}\in \tilde D_{N-1}^{n}} |f({\bf k})|^2=
\left( \frac N2\right)^n \sum _{{\bf r}\in \tilde D_{N-1}^{n}} h_{\bf r}
|a_{\bf r}|^2.
\]

{\bf AMDCT-4.}  This transform is given by the kernel
 \begin{equation}\label{AMDCT-4-1}
{\rm cos}^{4,-}_{\bf r}({\bf k})=|S_n|^{-1/2} \det \left( \cos
\frac{\pi (r_i+\frac12)(k_j+\frac12)}N\right)_{i,j=1}^n,
 \end{equation}
where ${\bf r}\in \tilde D_{N-1}^{n}$ and $k_j\in D_{N-1}$.
The orthogonality relation for these kernels is
 \begin{equation}\label{AMDCT-4-2}
\langle {\rm cos}^{4,-}_{\bf r}({\bf k}), {\rm cos}^{4,-}_{{\bf
r}'}({\bf k}) \rangle = |S_n|\sum_{{\bf k}\in \tilde D_{N-1}^{n}}
{\rm cos}^{4,-}_{\bf r}({\bf k}) {\rm cos}^{4,-}_{{\bf r}'}({\bf k})
 = \left(\frac N2\right)^n \delta_{{\bf r}{\bf r}'}.
 \end{equation}

This transform is given by the formula
 \begin{equation}\label{AMDCT-4-3}
f( {\bf k}){=}\sum _{{\bf r}\in \tilde D_{N-1}^{n}} a_{\bf r} {\rm
cos}^{4,-}_{\bf r}({\bf k}) ,\quad {\rm where}\quad a_{\bf
r}{=}|S_n| \left( \frac2N\right)^n \sum_{{\bf k}\in \tilde
D_{N-1}^{n}}  f({\bf k}) {\rm cos}^{4,-}_{\bf r}({\bf k}).
 \end{equation}
The Plancherel formula is
\[
 |S_n| \sum_{{\bf k}\in \tilde D_{N-1}^{n}} |f({\bf k})|^2=
\left( \frac N2\right)^n \sum _{{\bf r}\in \tilde D_{N-1}^{n}} |a_{\bf
r}|^2.
\]

\section{Other symmetric multivariate finite cosine transforms}
To each of the finite cosine transforms DCT-1, DCT-2, DCT-3, DCT-4
there corresponds a symmetric multivariate finite cosine
transform. We denote the corresponding transforms as SMDCT-1,
SMDCT-2, SMDCT-3, SMDCT-4. Below we give these transforms without
proof (proofs are the same as in the case of symmetric multivariate
finite cosine transforms of section 11). We fix a positive
integer $N$ and use the notation $\breve D_{N}^{n}$ for the subset of the
set $D_N^n\equiv D_N\times D_N\times \dots \times D_N$ ($n$ times)
with $D_N=\{ 0,1,2,\dots,N\}$ consisting of points ${\bf
r}=(r_1,r_2,\dots,r_n)$, $r_i\in {\mathbb Z}^{\ge 0}$, such that
\[
 N\ge r_1\ge r_2\ge \cdots \ge r_n\ge 0.
\]
The set $\breve D_{N}^{n}$ is an extension of the set
$\tilde D_{N}^{n}$ of
the previous section by adding points which are invariant with
respect of some elements of the permutation group $S_n$.

The set $D_N^n$ is obtained by action by elements of the group $S_n$
upon $\breve D_{N}^{n}$, that is, $D_N^n$ coincides with the set $\{
w\breve D_{N}^{n}; w\in S_n\}$. However, in
$\{ w\breve D_{N}^{n}; w\in
S_n\}$, some points are met several times. Namely, a point  ${\bf
k}_0\in \breve D_{N}^{n}$ is met $|S_{{\bf k}_0}|$ times in the set $\{ w
\breve D_{N}^{n}; w\in S_n\}$, where $|S_{{\bf k}_0}|$ is an order of the
subgroup $S_{{\bf k}_0}\subset S_n$ consisting of elements $w\in
S_n$ leaving ${\bf k}_0$ invariant.
\medskip

{\bf SMDCT-1.}  This transform is given by the kernel
 \begin{equation}\label{SMDCT-1-1}
{\rm cos}^{+}_{\bf r}({\bf k})\equiv {\rm cos}^{1,+}_{\bf r}({\bf
k})= |S_n|^{-1/2} {\det}^+\left( \cos \frac{\pi
r_ik_j}N\right)_{i,j=1}^n,
 \end{equation}
where  ${\bf k}=(k_1,k_2,\dots,k_n)$,
$k_i\in \{ 0,1,2,\dots,N\}$. The orthogonality relation for these kernels is
\[
\langle {\rm cos}^{+}_{\bf r}({\bf k}),{\rm cos}^{+}_{{\bf r}'}({\bf
k}) \rangle =|S_n|\sum_{{\bf k}\in \breve D_{N}^{n}} |S_{\bf
k}|^{-1} c_{\bf k} {\rm cos}^{+}_{\bf r}({\bf k}){\rm cos}^{+}_{{\bf
r}'}({\bf k})
\]
 \begin{equation}\label{SMDCT-1-2}
 =h_{\bf r}
\left(\frac N2\right)^n |S_{\bf r}|  \delta_{{\bf r}{\bf r}'},
 \end{equation}
where $S_{\bf r}$ is the subgroup of $S_n$ consisting of elements
leaving ${\bf r}$ invariant,
\[
 c_{\bf k}=c_1c_2\cdots c_n,\ \ \  h_{\bf k}=h_1h_2\cdots h_n,
\]
and $c_i$ and $h_j$ are such as in formula \eqref{DCT-1-3}.

This transform is given by the formula
 \begin{equation}\label{SMDCT-1-3}
f( {\bf k}){=}\sum _{{\bf r}\in \breve D_{N}^{n}} a_{\bf r}
 {\rm cos}^{+}_{\bf r}({\bf k}) ,
 \end{equation}
 where
\[
a_{\bf r}{=}\frac{|S_n|}{h_{\bf r}|S_{\bf r}|} \left(
\frac2N\right)^n \sum_{{\bf k}\in \breve D_{N}^{n}} |S_{\bf k}|^{-1}
 c_{\bf k} f({\bf k}) {\rm cos}^{+}_{\bf r}({\bf k}).
\]
The Plancherel formula is
\[
 |S_n|\sum_{{\bf k}\in \breve D_{N}^{n}} |S_{\bf k}|^{-1}
c_{\bf k} |f( {\bf k})|^2=
 \left(\frac N2\right)^n
  \sum _{{\bf r}\in \breve D_{N}^{n}} h_{\bf r} |S_{\bf r}| |a_{\bf r}|^2.
\]

This transform is in fact a variation of the symmetric multivariate
discrete cosine transforms from section 11.
\medskip

{\bf SMDCT-2.}  This transform is given by the kernel
 \begin{equation}\label{SMDCT-2-1}
 {\rm cos}^{2,+}_{\bf r}({\bf k})= |S_n|^{-1/2} {\det}^+\left(
\cos \frac{\pi (r_i+\frac12)k_j}N\right)_{i,j=1}^n,\ \ \
 {\bf r}\in \breve D_{N-1}^{n},
 \end{equation}
where $\breve D_{N-1}^{n}$ is the set $\breve D_{N}^{n}$ with $N$ replaced by
$N-1$ and ${\bf k}=(k_1,k_2,\dots, k_n)$, $k_i\in \{ 0,1,2,\dots, N-1\}$. The
orthogonality relation for these kernels is
\[
\langle {\rm cos}^{2,+}_{\bf r}({\bf k}),{\rm cos}^{2,+}_{{\bf
r}'}({\bf k}) \rangle =|S_n|\sum_{{\bf k}\in \breve D_{N-1}^{n}}
|S_{\bf k}|^{-1}c_{\bf k} {\rm cos}^{2,+}_{\bf r}({\bf k}) {\rm
cos}^{2,+}_{{\bf r}'}({\bf k})
\]
 \begin{equation}\label{SMDCT-2-2}
=\left(\frac N2\right)^n |S_{\bf r}| \delta_{{\bf r}{\bf r}'},
 \end{equation}
where  $c_{\bf k}=c_1c_2\cdots c_n$ and $c_j$ are such as in
\eqref{DCT-2-3}.

This transform is given by the formula
 \begin{equation}\label{SMDCT-2-3}
f( {\bf k})=\sum _{{\bf r}\in \breve D_{N-1}^{n}} a_{\bf r}
 {\rm cos}^{2,+}_{\bf r}({\bf k}) ,
 \end{equation}
 where
\[
a_{\bf r}=\frac{|S_n|}{|S_{\bf r}|} \left( \frac2N\right)^n
\sum_{{\bf k}\in \breve D_{N-1}^{n}} |S_{\bf k}|^{-1} c_{\bf k}
f({\bf k}) {\rm cos}^{2,+}_{\bf r}({\bf k}).
\]
The Plancherel formula is of the form
\[
 |S_n|\sum_{{\bf k}\in \breve D_{N-1}^{n}}
|S_{\bf k}|^{-1} c_{\bf k} |f( {\bf k})|^2=
 \left(\frac N2\right)^n
  \sum _{{\bf r}\in \breve D_{N-1}^{n}}|S_{\bf r}| |a_{\bf r}|^2.
\]

{\bf SMDCT-3.}  This transform is given by the kernel
  \begin{equation}\label{SMDCT-3-1}
{\rm cos}^{3,+}_{\bf r}({\bf k})= |S_n|^{-1/2} {\det}^+\left( \cos
\frac{\pi r_i(k_j+\frac12)}N\right)_{i,j=1}^n,
 \end{equation}
where ${\bf r}\in \breve D_{N-1}^{n}$.
The orthogonality relation for these kernels is
\[
\langle {\rm cos}^{3,+}_{\bf r}({\bf k}),{\rm cos}^{3,+}_{{\bf
r}'}({\bf k}) \rangle =|S_n|\sum_{{\bf k}\in \breve D_{N-1}^{n}}
|S_{\bf k}|^{-1} {\rm cos}^{3,+}_{\bf r}({\bf k}) {\rm
cos}^{3,+}_{{\bf r}'}({\bf k})
\]
\begin{equation}\label{SMDCT-3-2}
=h_{\bf
r} \left(\frac N2\right)^n |S_{\bf r}| \delta_{{\bf r}{\bf r}'},
 \end{equation}
where $h_{\bf r}=h_1h_2\cdots h_n$ and $h_i$ are such as in formula
\eqref{DCT-3-3}.

This transform is given by the formula
 \begin{equation}\label{SMDCT-3-3}
f( {\bf k})=\sum _{{\bf r}\in \breve D_{N-1}^{n}} a_{\bf r}
 {\rm cos}^{3,+}_{\bf r}({\bf k}) ,
 \end{equation}
where
\[
a_{\bf r}=\frac{|S_n|}{h_{\bf r}|S_{\bf r}|} \left( \frac2N\right)^n
\sum_{{\bf k}\in \breve D_{N-1}^{n}} |S_{\bf k}|^{-1} f({\bf k})
{\rm cos}^{3,+}_{\bf r}({\bf k}).
\]
The Plancherel formula is of the form
\[
 |S_n|\sum_{{\bf k}\in \breve D_{N-1}^{n}} |S_{\bf k}|^{-1} |f( {\bf k})|^2=
 \left(\frac N2\right)^n
  \sum _{{\bf r}\in \breve D_{N-1}^{n}}h_{\bf r} |S_{\bf r}|  |a_{\bf r}|^2.
\]

{\bf SMDCT-4.}  This transform is given by the kernel
 \begin{equation}\label{SMDCT-4-1}
{\rm cos}^{4,+}_{\bf r}({\bf k})= |S_n|^{-1/2} {\det}^+\left( \cos
\frac{\pi (r_i+\frac12)(k_j+\frac12)}N\right)_{i,j=1}^n,
 \end{equation}
where ${\bf r}\in \breve D_{N-1}^{n}$.
The orthogonality relation for these kernels is
\[
\langle {\rm cos}^{4,+}_{\bf r}({\bf k}), {\rm cos}^{4,+}_{{\bf
r}'}({\bf k}) \rangle =|S_n|\sum_{{\bf k}\in \breve D_{N-1}^{n}}
|S_{\bf k}|^{-1} {\rm cos}^{4,+}_{\bf r}({\bf k}) {\rm
cos}^{4,+}_{{\bf r}'}({\bf k})
 \]
\begin{equation}\label{SMDCT-4-2}
 = \left(\frac
N2\right)^n |S_{\bf r}| \delta_{{\bf r}{\bf r}'}.
 \end{equation}

This transform is given by the formula
 \begin{equation}\label{SMDCT-4-3}
f( {\bf k})=\sum _{{\bf r}\in \breve D_{N-1}^{n}} a_{\bf r}
 {\rm cos}^{4,+}_{\bf r}({\bf k}) ,
 \end{equation}
where
\[
a_{\bf r}= \left( \frac2N\right)^n  \frac{|S_n|}{|S_{\bf r}|}
\sum_{{\bf k}\in \breve D_{N-1}^{n}} |S_{\bf k}|^{-1} f({\bf k})
{\rm cos}^{4,+}_{\bf r}({\bf k}).
\]
The Plancherel formula is
\[
 |S_n|\sum_{{\bf k}\in \breve D_{N-1}^{n}} |S_{\bf k}|^{-1} |f( {\bf k})|^2=
 \left(\frac N2\right)^n
  \sum _{{\bf r}\in \breve D_{N-1}^{n}}|S_{\bf r}| |a_{\bf r}|^2.
\]

\section*{Acknowledgements}

We are grateful for the hospitality extended to A.K. at the
Center de Recherches Math\'ematiques, Universit\'e de Montr\'eal,
during the preparation of this paper.
His research was partially supported by Grant 14.01/016 of
the State Foundation of Fundamental Research of Ukraine.
We acknowledge also partial support for this work from
the National Science and Engineering Research Council of
Canada, MITACS, the MIND Institute of Costa Mesa, California, and
Lockheed Martin, Canada.

\end{document}